\documentclass[11pt]{article}
\usepackage{amsmath,amsthm,amssymb}
\usepackage{thmtools,thm-restate}
\usepackage[noadjust]{cite}
\usepackage{graphicx}
\usepackage{mathrsfs}
\usepackage{hyperref}
\usepackage{mathtools}
\hypersetup{colorlinks=true,linkcolor=blue,filecolor=blue,citecolor=blue}
\usepackage{cleveref}
\usepackage{makecell}
\usepackage{array}
\usepackage[defaultlines=3,all]{nowidow}
\usepackage{microtype}
\usepackage{multirow}
\usepackage{mdwtab}
\usepackage{comment}
\newtheorem{theorem}{Theorem}[section]
\newtheorem{lemma}[theorem]{Lemma}

\newtheorem{conjecture}[theorem]{Conjecture}
\newtheorem*{conjecture*}{Conjecture}

\newtheorem{claim}{Claim}
\newtheorem{observation}{Observation}

\begin{document}

\title{\bf  Semistrong edge coloring and $(0,1)$-relaxed strong edge coloring of graphs
}
\author{Yuquan Lin and Wensong Lin\footnote{Corresponding author. E-mail address: wslin@seu.edu.cn}\\
{\small School of Mathematics, Southeast University, Nanjing 210096, P.R. China}}
\date{}
\maketitle

\vspace*{-1cm} \setlength\baselineskip{6.0mm}
\bigskip

\begin{abstract}

In this work, we study two relaxations of the well-known strong edge coloring.
A \emph{semistrong edge coloring} of a graph $G$ is an edge coloring  in which every color class forms a matching $M$ such that every edge of $M$ is incident with (at least) one vertex  of degree 1 in  the subgraph of $G$ induced by the  vertices covered by $M$.
For any two nonnegative integers $s$ and $t$,  an \emph{$(s,t)$-relaxed  strong edge coloring}   of $G$  is  an edge coloring in which,  for  every edge  $e$ of $G$, at most $s$ edges at distance 1  and at most $t$ edges at distance 2 from $e$   receive the same color as $e$. 
 The corresponding chromatic indices are defined accordingly.
 
We  confirm a recent conjecture of  Lu{\v{z}}ar, Mockov{\v{c}}iakov{\'a}, and  Sot{\'a}k [J. Graph Theory 105 (2024) 612–632],
 which asserts that 
every connected graph $G$ with  maximum degree $\Delta$ ($\ge3$), except for $K_{\Delta,\Delta}$, has a semistrong chromatic index at most $\Delta^2-1$.
This is achieved by constructing an edge coloring of $G$ using at most $\Delta^2-1$ colors that is simultaneously  semistrong  and  $(0,1)$-relaxed strong.
Consequently, every such graph  also has $(0,1)$-relaxed strong chromatic index at most $\Delta^2-1$.
\end{abstract}



\section{Introduction}
\label{sec:1}
In this paper, we only consider  finite undirected simple connected graphs.  Let $G=(V(G), E(G))$ be a  graph.
For $v\in V(G)$, let $N(v)=\{u\in V(G): uv\in E(G)\}$ denote the open neighborhood of $v$ and  $d(v)=|N(v)|$ be the degree of $v$. Let  $\Delta=\max\limits_{v\in V(G)}d(v)$ denote the maximum degree of $G$.  
For  $M\subseteq E(G)$, we denote by $G_M$  
 the subgraph of $G$ induced by the   endvertices  of the edges in $M$.

Given two positive integers $i$ and $j$, denote by  $C_{i}$  the cycle on  $i$ vertices,  by  $K_{i}$  the complete graph on $i$ vertices, and
by $K_{i,j}$ the complete bipartite graph with parts of sizes $i$ and $j$, respectively.
For convenience, we use the abbreviation $[1,i]$ for $\{1,2, \dots, i\}$. 

Let $e$ and $e'$ be two edges of $G$. If $e$ and $e'$ are adjacent to each other, we say that the distance between $e$ and $e'$ is $1$, and if they are not adjacent but both of them are adjacent to a common edge, we say they are at distance $2$.
 An {\em induced matching} (also called a {\em strong matching})  $M$ of $G$ is a matching such that no two edges of $M$ are at distance $1$ or $2$ in $G$. In other words, a matching $M$ of $G$ is induced if each vertex in  $G_M$ is of degree 1.


  Given a positive integer $k$,  a {\em strong $k$-edge-coloring} of $G$ is  an assignment of $k$ colors to the edges of $G$ such that every color class forms  an induced matching.
 The {\em strong chromatic index} $\chi'_{s}(G)$ of $G$  is the minimum integer $k$ for which $G$ admits a strong $k$-edge-coloring.
 
The concept of strong edge coloring, first introduced by Fouquet and Jolivet \cite{FJ1983}, can be used to model the conflict-free channel assignment problem in radio networks \cite{R1997,NKGB2000}. In 1985,   Erd\H{o}s and  Ne\v{s}et\v{r}il \cite{E1988,EN1989} proposed the following conjecture on the upper bound of $\chi'_{s}(G)$ in terms of the maximum degree $\Delta$. If true, the conjectured bound is the best possible.

\begin{conjecture}\label{Conj-EN}
	$($Erd\H{o}s and Ne\v{s}et\v{r}il \cite{E1988,EN1989}$)$
	If $G$ is a graph with maximum degree $\Delta$, then
	\begin{equation*}
	\chi'_{s}(G)\le\begin{cases}
			\begin{array}{cl}
				\dfrac{5}{4}\Delta^{2},                                 & \text{if}\  \Delta \ \text{is even,} \\
				\dfrac{5}{4}\Delta^{2}-\dfrac{1}{2}\Delta+\dfrac{1}{4}, & \text{if}\  \Delta\  \text{is odd.}
			\end{array}
	\end{cases}	
	\end{equation*}
\end{conjecture}

This conjecture is considered one of the most important  in the study of strong edge coloring. 
Over the past few decades, numerous studies on strong edge coloring have been motivated by this conjecture; however, little progress has been made toward a direct proof.
Only the case $\Delta\le 3$ was confirmed completely   by Andersen \cite{A1992} in 1992, and independently by Hor\'{a}k, Qing, and Trotter \cite{HQT1993} in 1993. 
Beyond this, the problem remains widely open.

 For sufficiently large $\Delta$, Molloy and Reed \cite{MR1997} first proved in 1997 that $\chi'_{s}(G)\le 1.998\Delta^{2}$  using  probabilistic techniques.
 This bound was improved  to   $1.93\Delta^{2}$  by Bruhn and Joos \cite{BJ2015}  in 2015, and  further strengthened  to $1.835\Delta^{2}$  by Bonamy, Perrett and Postle \cite{BPP2022} in 2022.   
 The current best known upper bound  $1.772\Delta^{2}$ was established  by Hurley, de Joannis de Verclos and Kang \cite{HdK2021} in 2021. 
The results 
 mentioned above apply a similar proof method, but this method has its  limitations, so  the  best possible  coefficient by far  is still not very close to the   objective of 1.25.


 It seems difficult to  prove 
Conjecture \ref{Conj-EN} directly. 
Recently, a lot of attention has been paid to various  variants of  strong edge coloring (see, e.g., \cite{BR2018,GHHL2005,GH2005,HL2017,HLL2022}).
In 2005, Gy{\'a}rf{\'a}s and  Hubenko \cite{GH2005} introduced the  {\em semistrong edge coloring} by relaxing the strong  (induced)  matching requirement of strong edge coloring to the notion of the \emph{semistrong matching}, that is, a matching $M$ in which every edge is incident with (at least) one vertex  of degree 1 in the induced subgraph $G_M$.
Formally,  a \emph{semistrong $k$-edge-coloring}  of $G$ is an edge coloring using at most $k$ colors  in which every color class is a semistrong matching. The minimum integer $k$ such that $G$ has a semistrong $k$-edge-coloring is called the {\em semistrong chromatic index} of $G$, denoted by $\chi'_{ss}(G)$.
It is obvious that   $\chi'_{ss}(G)\le \chi'_{s}(G)$. Gy{\'a}rf{\'a}s and  Hubenko \cite{GH2005}  proved that  if $G$ is a {\em Kneser graph} or a {\em subset graph}, then $\chi'_{ss}(G)=\chi'_{s}(G)$.

Recently,  Lu{\v{z}}ar, Mockov{\v{c}}iakov{\'a}, and  Sot{\'a}k  \cite{LMS2022} revived  the semistrong edge coloring and further investigated its properties. 
They showed that complete graphs and  complete bipartite graphs are two additional classes of graphs for which the semistrong and strong chromatic indices coincide.
(To be precise, $\chi'_{ss}(K_{n})=\chi'_{s}(K_{n})=\binom{n}{2}$ and $\chi'_{ss}(K_{m,n})=\chi'_{s}(K_{m,n})=mn$.)
Moreover, they revealed the fact that, according to the work of Diwan \cite{D2019} and  the work of  Faudree,   Schelp,   Gy{\'a}rf{\'a}s,  and   Tuza \cite{FGST1990},  it can be concluded that
 $\chi'_{ss}(Q^{n})=\chi'_{s}(Q^{n})=2n$ for any {\em $n$-dimensional cube} $Q^{n}$ with $n\ge2$.

\begin{figure}[htbp]  
 	\centering	
 	\begin{minipage}[t]{7cm}
 		\centering
 		\resizebox{2.9cm}{2.5cm}{\includegraphics{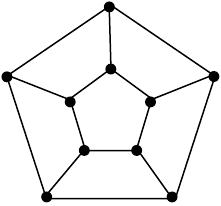}}
 		\caption{The graph $5$-prism.}
 		\label{fig:5prism}
 	\end{minipage}
 	\begin{minipage}[t]{7cm}
 		\centering
	 		\resizebox{3.6cm}{2.3cm}{\includegraphics{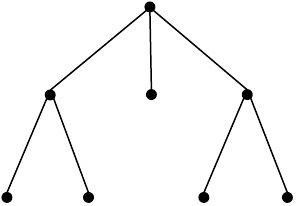}}
 		\caption{The graph $T_{0}$.}
 		\label{fig:T0}
 	\end{minipage}
 \end{figure}

In  \cite{LMS2022}, the authors also proved that 
 $\chi'_{ss}(G)\le \Delta^{2}$ for every graph $G$ with maximum degree $\Delta$. Moreover, for the case $\Delta=3$,
 they improved the bound 9 to 8 for every connected graph $G$ that is not isomorphic to $K_{3,3}$, where the $5$-prism (as shown in Figure \ref{fig:5prism}) shows
 the sharpness of the upper bound 8. 
At the end of their paper, they proposed the  following conjecture.

\begin{conjecture}\label{conjecture-semi} 
	$($Lu{\v{z}}ar, Mockov{\v{c}}iakov{\'a},  Sot{\'a}k  \cite{LMS2022}$)$
		For every connected graph $G$ with maximum degree $\Delta$, distinct from  $K_{\Delta,\Delta}$, it holds that $\chi'_{ss}(G)\le \Delta^{2}-1$. 		
\end{conjecture}

This paper settles this conjecture by proving the following two theorems.

\begin{theorem}\label{main-degree2}
	Let $G$ be a connected graph with maximum degree $2$. If $G$ is not isomorphic to   $C_{4}$ or $C_{7}$, then $\chi'_{ss}(G)\le 3$; otherwise,  $\chi'_{ss}(G)= 4$. 	
\end{theorem}

\begin{theorem}\label{main-degree3+}
Let $G$ be a  connected graph with maximum degree $\Delta \ge 3$. If  $G$ is not isomorphic to  $K_{\Delta,\Delta}$, then $\chi'_{ss}(G)\le \Delta^{2}-1$; otherwise, $\chi'_{ss}(G)=\Delta^{2}$.
\end{theorem}



It should be pointed out that  different relaxations of strong edge coloring   may be related to each other.
For example,  the  $(s,t)$-relaxed strong edge coloring,  which was first proposed by  He and Lin  \cite{HL2017} in 2017, is suitable for the channel assignment problem with  limited channel resources in wireless radio networks.
For any nonnegative integers $s$, $t$, and $k$, an  \emph{$(s,t)$-relaxed strong $k$-edge-coloring}  of $G$ is an assignment of $k$ colors to edges of $G$, such that for   every  edge  $e$ of $G$, at most $s$ edges at distance $1$  and at most $t$ edges at distance $2$ from $e$   receive the same color as $e$. 
The \emph{$(s,t)$-relaxed  strong chromatic index} $\chi'_{(s,t)}(G)$  of $G$ is the minimum integer $k$ such that $G$ admits an  $(s,t)$-relaxed strong $k$-edge-coloring.

In  \cite{HL2017}, He and Lin studied   $(s,t)$-relaxed strong edge coloring of {\em trees} and constructed a $(0,\Delta-1)$-relaxed strong $(\Delta+1)$-edge-coloring for any given tree $T$ with maximum degree $\Delta$.
Later, Lu{\v{z}}ar, Mockov{\v{c}}iakov{\'a} and  Sot{\'a}k  \cite{LMS2022} observed that the coloring provided by He and Lin is also a semistrong edge coloring, which immediately implies  $\chi'_{ss}(T)\le \Delta+1$ for any tree $T$.
Moreover, they proved in  \cite{LMS2022} that for any graph $G$, there exists an edge coloring using at most $\Delta^{2}$ colors  that is both a semistrong edge coloring and a $(0,1)$-relaxed strong edge coloring. 
(Hence, every graph $G$ with maximum degree $\Delta$ satisfies $\chi'_{(0,1)}(G)\le \Delta^{2}$.)

Inspired by their work in \cite{LMS2022}, in addressing \Cref{conjecture-semi}, we construct an edge coloring that is simultaneously a semistrong edge coloring and a $(0,1)$-relaxed strong edge coloring, thereby establishing the following result.

\begin{theorem}\label{main-relaxed}
Let $G$ be a  connected graph with maximum degree $\Delta \ge 2$. If  $G$ is not isomorphic to $C_{7}$,  then $\chi'_{(0,1)}(G)\le \Delta^{2}-1$; otherwise,  $\chi'_{(0,1)}(G)=4$.
\end{theorem}

\noindent{\bf Remark 1.} The semistrong chromatic  index and the $(s,t)$-relaxed  strong chromatic index of a graph $G$ are not comparable.
For instance, for the cycle $C_{4}$, it holds that $\chi'_{ss}(C_{4})=4> \chi'_{(0,1)}(C_{4})=2$, whereas for the cycle
 $C_{7}$,  there is $\chi'_{ss}(C_{7})= \chi'_{(0,1)}(C_{7})=4$. 
While for the graph $T_{0}$ shown  in \Cref{fig:T0}, we have  $\chi'_{ss}(T_{0})=3< \chi'_{(0,1)}(T_{0})=4$.\\

\noindent{\bf Remark 2.}  For the strong chromatic index of a graph, the upper bound in  \Cref{Conj-EN} is $1.25\Delta^{2}$, while the best known result for  large $\Delta$ is  $1.772\Delta^{2}$   provided by Hurley, de Joannis de Verclos and Kang \cite{HdK2021}. 
In contrast, our bounds for both semistrong chromatic index and $(0,1)$-relaxed strong chromatic index are $\Delta^{2}-1$.
This implies that, even a slight relaxation  can lead to a substantial reduction in the number of colors required.

\medskip
The remainder of this paper is dedicated to the proofs of  \Cref{main-degree2,main-degree3+,main-relaxed}. It  is organized as follows.
In  \Cref{sec:2}, we deal with the case that the maximum degree of  $G$ is $2$  and the case that $G$ is isomorphic to  $K_{\Delta,\Delta}$. 
In the next two sections, we assume that $G$ is a connected graph with maximum degree $\Delta \ge 3$ and not isomorphic to $K_{\Delta,\Delta}$.
In \Cref{sec:3}, after stating necessary definitions and  notation, we introduce some auxiliary results that will play a crucial role in the subsequent proofs.
\Cref{sec:4} is the main body of the proofs, where we construct  an edge coloring of $G$ using at most $\Delta^2-1$ colors that is both  semistrong  and  $(0,1)$-relaxed strong.
 Finally, we  summarize our results and suggest some future research directions in  \Cref{sec:5}.

\section{The proofs of two special cases}
\label{sec:2}

Let $G$ be a connected graph with maximum degree $\Delta\ge2$.
In this section, we  consider the semistrong chromatic index and the $(0,1)$-relaxed strong chromatic index of $G$ when  $\Delta=2$ (see \Cref{lemma:delta2}) and when $G$ is isomorphic to  $K_{\Delta,\Delta}$ (see Lemma~\ref{lemma:Kn,n}), respectively.

We first look at the case where $\Delta=2$.


\begin{lemma}\label{lemma:delta2}
	Let $G$ be a connected graph with maximum degree $2$.  Then,
    \begin{itemize}
        \item $\chi'_{ss}(C_{4})=4$ and  $\chi'_{(0,1)}(C_{4})=2$; 
        \item   $\chi'_{ss}(C_{7})=4$ and  $\chi'_{(0,1)}(C_{7})=4$;
        \item  if $G$ is not  isomorphic to  $C_{4}$ or $C_{7}$, then $\chi'_{ss}(G)\le 3$ and  $\chi'_{(0,1)}(G)\le 3$.
    \end{itemize}
\end{lemma}
\begin{proof}
Note that all edges of  $C_{4}$ must receive different colors in any semistrong edge coloring, hence  $\chi'_{ss}(C_{4})=4$.  It is obvious that $\chi'_{(0,1)}(C_{4})=2$.

Observe that any semistrong matching in $C_{7}$  consists of at most two edges,
we have $\chi'_{ss}(C_{7})\ge 4$.
Since  a  semistrong edge coloring of $C_{7}$ using $4$ colors is easy to get, it holds that $\chi'_{ss}(C_{7})\le 4$.
Hence,  $\chi'_{ss}(C_{7})=4$. 
 Similarly, it is straightforward to check that $\chi'_{(0,1)}(C_{7})=4$.

Assume that  $G$ is not  isomorphic to  $C_{4}$ or $C_{7}$. Note that $G$ is either a path or a cycle.	
Suppose that $G$ is a path with $n$ vertices.
Without loss of generality, label the vertices of $G$ as $v_{1},v_{2},\ldots ,v_{n}$, and let $e_{i}=v_{i}v_{i+1}$ for each $i\in[1,n-1]$. 
Define an edge coloring $\phi$ of $G$ by setting  $\phi(e_{i})=i \bmod 3$ for each $i\in[1,n-1]$.
This gives an edge coloring  of $G$ using  at most 3 colors,  which is both semistrong and  $(0,1)$-relaxed strong.

Next, we suppose that $G$ is a cycle $C_{n}=v_{1}v_{2}\ldots v_{n}$ with  $n\ge3$ and $n\notin\{4,7\}$.
Denote the edge $v_{i}v_{i+1}$ by $e_{i}$ for each $i\in[1,n-1]$ and the edge $v_{n}v_{1}$ by $e_{n}$.
If $n \equiv 1(\bmod\ 3)$, define an edge coloring 
$\phi$ by setting $\phi(e_{i})=i \bmod 3$ for each $i\in[1,n-4]$, and  letting  $\phi(e_{n-3})=2$, $\phi(e_{n-2})=1$,  $\phi(e_{n-1})=0$, and $\phi(e_{n})=2$.
 Otherwise, let $\phi(e_{i})=i \bmod 3$ for each $i\in[1,n]$.
It is easy to check that, in both cases, we obtain a semistrong edge coloring  $\phi$ of $G$ using 3 colors which is also a $(0,1)$-relaxed strong edge coloring.  
 Therefore, Lemma \ref{lemma:delta2} is proved.
\end{proof}

 \Cref{main-degree2} follows directly from the above lemma.
We proceed to analyze the case where $G$ is isomorphic to  $K_{\Delta,\Delta}$.

\begin{lemma}\label{lemma:Kn,n}
	 $\chi'_{ss}(K_{\Delta,\Delta})= \Delta^{2}$ and  $\chi'_{(0,1)}(K_{\Delta,\Delta})= \lceil \frac{\Delta^{2}}{2}\rceil$. 
\end{lemma}
\begin{proof}
Recall that all edges of  $C_{4}$ must receive different colors in any semistrong edge coloring, and thus  all edges in $K_{\Delta,\Delta}$ must receive different colors in any semistrong edge coloring.
Hence, $\chi'_{ss}(K_{\Delta,\Delta})= \Delta^{2}$.

Since any two edges of $K_{\Delta,\Delta}$ are at distance 1 or 2, every color class in a $(0,1)$-relaxed strong edge coloring of $K_{\Delta,\Delta}$ contains at most two edges. Thus, we have  $\chi'_{(0,1)}(K_{\Delta,\Delta})\ge \lceil \frac{\Delta^{2}}{2}\rceil$. 
Denote the two parts of $K_{\Delta,\Delta}$ by $U=\{u_{1},u_{2},\ldots,u_{\Delta}\}$ and $V=\{v_{1},v_{2},\ldots,v_{\Delta}\}$, respectively.
We define an edge coloring  $\phi$ of $G$ as follows:
let $\phi(u_{i}v_{j})=\phi(u_{j}v_{i})=\alpha_{i,j}$   for any two different integers $i,j\in[1,\Delta]$, 
and let  $\phi(u_{i}v_{i})=\beta_{\lceil \frac{i}{2}\rceil}$ for each $i\in[1,\Delta]$.
It is clear that $\phi$ is 
a  $(0,1)$-relaxed strong edge coloring using $\tbinom{\Delta}{2}+\lceil \frac{\Delta}{2}\rceil=\lceil \frac{\Delta^{2}}{2}\rceil$ colors,  and so $\chi'_{(0,1)}(K_{\Delta,\Delta})\le \lceil \frac{\Delta^{2}}{2}\rceil$. 
Hence, 
$\chi'_{(0,1)}(K_{\Delta,\Delta})= \lceil \frac{\Delta^{2}}{2}\rceil$. 
The lemma holds.
\end{proof}

By Lemmas~\ref{lemma:delta2} and~\ref{lemma:Kn,n}, it suffices to  complete the proofs of   \Cref{main-degree3+,main-relaxed} by  proving the following theorem.

\begin{theorem}\label{lemma:delta3+}
Let $G$ be a connected  graph with maximum degree $\Delta \ge 3$. If $G$ is not isomorphic to  $K_{\Delta,\Delta}$, then $\chi'_{ss}(G)\le \Delta^{2}-1$ and $\chi'_{(0,1)}(G)\le \Delta^{2}-1$. 
\end{theorem}

The remainder of this paper is devoted to the proof of \Cref{lemma:delta3+}.

\section{Preliminaries and notation}
\label{sec:3}
In this section, we  introduce some notation and preliminary facts that will be used in our proofs.
We usually use $\alpha$, $\beta$, $\gamma$ to denote colors, and $\phi$, $\psi$, $\sigma$ to denote edge colorings.  Sometimes,  we simply write “coloring” instead of “edge coloring”.

Let $G$ be a graph and let $\phi$ be an edge   coloring  of $G$.
  For $S\subseteq E(G)$, we denote by $\phi(S)$ the set of colors assigned to the edges in $S$ under  $\phi$.

For any two edges $e,f\in E(G)$,
 we say that  $f$ is a {\em $1$-neighbor} (resp.  {\em $2$-neighbor}) of $e$ if  $f$ and $e$ are at distance 1 (resp.  2),  and  $f$ is a {\em $2^{-}$-neighbor}  of $e$ if  they  are at distance 1 or 2. 
For any $e\in E(G)$,  we use $C_{e}^{\Delta}$ to denote the set of $1$-neighbors of $e$  lying on a common $3$-cycle with $e$.

For each edge $e=uv\in E(G)$, we denote by $N(e)$ (resp. $N^{2}(e)$) the set of $1$-neighbors (resp.   $2$-neighbors) of  $e$, and by $N^{2-}(e)$ the set  of  $2^{-}$-neighbors of $e$.
It is obvious that $N(e)\cap N^{2}(e)=\emptyset$ and  $N^{2-}(e)=N(e)\cup N^{2}(e)$.
Let $N[e]=N(e)\cup \{e\}$ and  $N^{2-}[e]=N^{2-}(e)\cup \{e\}$.
Similarly, let $N_{u}(e)$ denote the set of $1$-neighbors of  $e$ having $u$ as an  endvertex and  $N_{u}^{2}(e)$ the set of  $2$-neighbors of $e$  being adjacent to some edge in $N_{u}(e)$. Denote by $N_{u}^{2-}(e)$ the set of  edges in $N_{u}(e)$ or $N_{u}^{2}(e)$.
Moreover,  let $N_{u}[e]=N_{u}(e)\cup \{e\}$ and  $N_{u}^{2-}[e]=N_{u}^{2-}(e)\cup \{e\}$.

\begin{figure}[htbp]  
		\centering
		\resizebox{11.2cm}{2.5cm}{\includegraphics{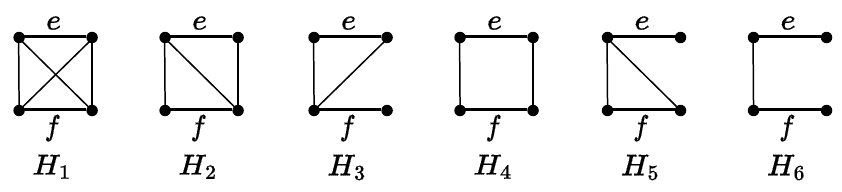}}
		\caption{The six possible configurations  of the induced subgraph $G_{\{e,f\}}$.}
		\label{fig:2-edge}
\end{figure}

Let $e$ and $f$ be two edges of $G$ at distance 2.
As shown in \Cref{fig:2-edge}, there are six possible configurations of the subgraph $G_{\{e,f\}}$ (induced by the four endvertices of $e$ and $f$).
If $G_{\{e,f\}}$ is the same as the graph $H_{i}$ shown in  \Cref{fig:2-edge}, then  we say that $f$ is   a {\em $2$-neighbor  of  Type $i$}  of $e$, where $i\in[1,6]$.
We denote by $T_{i}(e)$  the set of  $2$-neighbors  of  Type $i$ of $e$.  
It is clear that $N^{2}(e)=\cup_{i=1}^{6}T_{i}(e)$ and $T_{i}(e)\cap T_{j}(e)=\emptyset$ for any two different integers $i,j\in[1,6]$.
We further define
 $F(e)=N(e)\cup(\cup_{i=1}^{5}T_{i}(e))$.
 Hence, we have $N^{2-}(e)=F(e)\cup T_{6}(e)$.
Note that  $f\in F(e)$ if and only if $e\in F(f)$, and that $f\in  T_{6}(e)$ if and only if $e\in  T_{6}(f)$.

Based on the above definitions, we immediately observe the following.

\begin{observation}\label{Obser:3-cycle}
Let $e=uv$ be  an edge of $G$. If $C_{e}^{\Delta}=\emptyset$, then $T_{1}(e)=T_{2}(e)=T_{3}(e)=\emptyset$. 
\end{observation}

Note that  in any semistrong edge coloring of $G$,
every edge $e$ must receive a color different from that of every edge $f\in F(e)$. We now make another useful observation.

\begin{observation}\label{Obser:shangjie} 
	Let $G$ be a graph with maximum degree $\Delta$. Then, 	for any edge $e=uv$ of $G$,   
	$$|F(e)|\le \Delta^{2}-1-\frac{1}{2} |C_{e}^{\Delta}|-|T_{1}(e)|-\frac{1}{2}|T_{2}(e)|-\frac{1}{2}|T_{6}(e)|.$$	
Moreover, if  equality holds, then every vertex in $ N(u)\cup N(v)$ is of degree $\Delta$ in $G$.
\end{observation}
\begin{proof}
Let $e=uv$ be an edge of $G$. 
On the one hand,
according to the partition of its $2$-neighbors,
 it is not difficult to see that
		\allowdisplaybreaks[4] 
		\begin{align*}	
	&	\sum_{w\in N(u)\setminus\{v\} }(d(w)-1)+
		\sum_{w\in N(v)\setminus\{u\} }(d(w)-1)\\
		\ge & \	
	|C_{e}^{\Delta}|+	4|T_{1}(e)|+3|T_{2}(e)|+2|T_{3}(e)|+2|T_{4}(e)|+2|T_{5}(e)|+|T_{6}(e)|\\
		= & \ 	
		2|\cup_{i=1}^{5}T_{i}(e)|+ 	|C_{e}^{\Delta}|+	2|T_{1}(e)|+|T_{2}(e)|+|T_{6}(e)|.
	\end{align*}	
On the other hand,
since $G$ is a graph with maximum degree $\Delta$, we have
	\begin{align}	
\sum\limits_{w\in N(u)\setminus\{v\} }(d(w)-1)+
\sum\limits_{w\in N(v)\setminus\{u\} }(d(w)-1) \le& 	 \	2(\Delta-1)^{2}.
\end{align}
Combining the above two inequalities, it holds that
	\begin{align*}	
			|\cup_{i=1}^{5}T_{i}(e)|\le& 	 \	(\Delta-1)^{2}-\frac{1}{2} |C_{e}^{\Delta}|-|T_{1}(e)|-\frac{1}{2}|T_{2}(e)|-\frac{1}{2}|T_{6}(e)|.
\end{align*}
Recall that $F(e)=N(e)\cup(\cup_{i=1}^{5}T_{i}(e))$ and $|N(e)|\le 2(\Delta-1)$, it is easy to check that
 	\begin{align}	
 |F(e)| 
 \le &\ 
  \Delta^{2}-1-\frac{1}{2} |C_{e}^{\Delta}|-|T_{1}(e)|-\frac{1}{2}|T_{2}(e)|-\frac{1}{2}|T_{6}(e)|.
 \end{align}
It is obvious that   if (2) is an  equality, then 
(1) must be an  equality.
This implies that,  for every $w\in N(u)\cup N(v)$, we must have $d(w)=\Delta$.  This completes the proof of   \Cref{Obser:shangjie}.
\end{proof}

Let $p$ be a  positive integer.  We denote by $\mathcal{G}_{p}$  the family of $p$-regular graphs  $G$  with $2p$ vertices in which there exists  an edge $e=uv\in E(G)$ such that $N(u)\cup N(v)=V(G)$.
Note that  every graph in $\mathcal{G}_{p}$ is connected.
Note also that 
$K_{p,p}\in \mathcal{G}_{p}$, and that
$|E(G)|=p^{2}$ for every $G\in \mathcal{G}_{p}$.

We are now ready to establish an essential ingredient needed in the next section.

\begin{lemma}\label{lemma:Gp-iff}
	Let $G$ be a connected graph with maximum degree $\Delta$.
	Then, $G$ contains an edge $e$ with $|F(e)|= \Delta^{2}-1$ if and only if $G \in \mathcal{G}_{\Delta}$.
\end{lemma}
\begin{proof}
	Let $e=uv$ be an edge of  $G$ with $|F(e)|= \Delta^{2}-1$. 
		By \Cref{Obser:shangjie},
	$C_{e}^{\Delta}=T_{1}(e)=T_{2}(e)=T_{6}(e)=\emptyset$ and  $d(w)=\Delta$ for each $w\in N(u)\cup N(v)$.
		At this time, we must have $V(G)=N(u)\cup N(v)$, as otherwise  since $G$ is connected, 
	there exists a vertex $x\in V(G)\setminus (N(u)\cup N(v))$ being adjacent to some vertex $w\in (N(u)\cup N(v))\setminus\{u,v\}$. Then  $xw$ is a $2$-neighbor of Type 6 of $e$ and so obtain a contradiction.  Therefore,	 we have $|V(G)|=2\Delta$  and thus	
	$G \in \mathcal{G}_{\Delta}$.

Suppose that $G \in \mathcal{G}_{\Delta}$.
Let   $e=uv$ be an edge of $G$ such that $N(u)\cup N(v)=V(G)$.
Note that $|E(G)|=\Delta^{2}$ and $|N(e)|=2(\Delta-1)$.
Note also that any edge $f\in E(G)\setminus N[e]$ is a $2$-neighbor of $e$, and thus
 $|N(e)|+|N^{2}(e)|=\Delta^{2}-1$.
Since $N(u)\cup N(v)=V(G)$,  we have $T_{6}(e)=\emptyset$ and so $|N(e)|+|N^{2}(e)|=|N(e)|+|\cup_{i=1}^{5}T_{i}(e)|=|F(e)|= \Delta^{2}-1$. 
Hence, the lemma holds.
\end{proof}

\section{The proof of  \Cref{lemma:delta3+}}
\label{sec:4}
In this section, we prove
that if $G$ is a connected graph with maximum degree $\Delta\ge3$ that is not isomorphic to $K_{\Delta,\Delta}$, then $\chi'_{ss}(G)\le \Delta^{2}-1$ and $\chi'_{(0,1)}(G)\le \Delta^{2}-1$.
The proof is divided into two cases depending on whether such a graph $G$ belongs to $\mathcal{G}_{\Delta}$,  the family of $\Delta$-regular graphs  $G$  on $2\Delta$ vertices that contain  an edge $e=uv\in E(G)$ with $N(u)\cup N(v)=V(G)$. 
These two cases are treated in Lemmas~\ref{lemma:G in Gp} and~\ref{lemma:G notin Gp}, respectively.
In both cases, we construct an edge coloring of $G$ using at most $\Delta^{2}-1$ colors, which is
 simultaneously semistrong and $(0,1)$-relaxed strong.

Let $G$ be a connected graph with maximum degree $\Delta\ge3$ that is not isomorphic to $K_{\Delta,\Delta}$.
 We first look at the case $G\in \mathcal{G}_{\Delta}$.


\begin{lemma} \label{lemma:G in Gp}
	 Let $G$ be a  connected graph with maximum degree $\Delta\ge3$ 
     that is not isomorphic to $K_{\Delta,\Delta}$.
     If $G\in  \mathcal{G}_{\Delta}$, 
     then $\chi'_{ss}(G)\le \Delta^{2}-1$ and $\chi'_{(0,1)}(G)\le \Delta^{2}-1$. 
\end{lemma}
\begin{proof}
Let $e$ be an  edge of $G$ with $N(u)\cup N(v)=V(G)$. Because $G$ belongs to $\mathcal{G}_{\Delta}$ and  is not isomorphic to $K_{\Delta,\Delta}$, there exist two distinct vertices $u'\in N(u)\setminus \{v\}$ and $v'\in N(v)\setminus \{u\}$ such that $u'v'\notin E(G)$. This implies that $uu'$ and $vv'$ do not lie on a common $4$-cycle. 
Notice that $|E(G)|=\Delta^{2}$,
 a semistrong $(\Delta^{2}-1)$-edge-coloring of $G$ can be easily obtained by coloring the two edges  $uu'$ and $vv'$ with the same color 1 and the remaining $\Delta^{2}-2$ edges with the other $\Delta^{2}-2$ colors. This coloring is obviously also a $(0,1)$-relaxed strong edge coloring of $G$. 
Therefore, $\chi'_{ss}(G)\le \Delta^{2}-1$ and $\chi'_{(0,1)}(G)\le \Delta^{2}-1$.
\end{proof}

The remainder of this section is devoted to the case  $G\notin  \mathcal{G}_{\Delta}$. Suppose that  $G\notin  \mathcal{G}_{\Delta}$.
It follows from Observation \ref{Obser:shangjie} and Lemma \ref{lemma:Gp-iff} that  $|F(e)|\le \Delta^{2}-2$ for every $e\in E(G)$.
Hence, the greedy algorithm, coloring 
the edges of $G$ one by one in any  order, will produce an edge coloring with  at most $\Delta^{2}-1$  colors, in which every edge $e$ receives a color distinct from all colors of edges in $F(e)$.
We call such a coloring   \emph{good}.

Let $\phi$ be a  good coloring  of $G$.
For an edge  $e$ of $G$, if  it has at least two $2$-neighbors  with  the same color as it under $\phi$, then we call it a {\em bad edge} with respect to $\phi$. 
For a $2$-neighbor $f$ of $e$ with $\phi(e)=\phi(f)$, we call them a {\em bad pair}  with respect to $\phi$. 
We denote by $\kappa_{1}(\phi)$  (resp. $\kappa_{2}(\phi)$)  the number of the bad edges (resp. the bad pairs)  with respect to $\phi$ in $G$.
Similarly, we use  $\kappa_{1}(\phi,\alpha)$  (resp. $\kappa_{2}(\phi,\alpha)$)  to denote the number of the bad edges  (resp. the bad pairs)  being colored the color $\alpha$ with respect to $\phi$ in $G$.

Based on the above definitions, we immediately observe the following.
\begin{observation}\label{Obser:no bad edge}
	Let $\phi$ be a good  coloring of a graph $G$.  If  no edge of $G$ is a bad edge with respect to $\phi$, then $\phi$ is both  a semistrong edge coloring and 
	a $(0,1)$-relaxed strong edge coloring.
\end{observation}

Among all good colorings of $G$, we refer to a coloring  with the fewest bad edges  as the \emph{$1$-optimal} coloring of $G$. Moreover,  if a $1$-optimal coloring has the least number of bad pairs among all $1$-optimal colorings of $G$, then we call it a  \emph{$2$-optimal} coloring of $G$. 

\medskip
In what follows, we prove that any 
 $2$-optimal coloring $\phi$ of $G$ contains no bad edges 
 with respect to $\phi$ (that is, $\kappa_{1}(\phi)=0$). 
 Consequently,  by  \Cref{Obser:no bad edge},
$\phi$ is both  a semistrong edge coloring and 
a $(0,1)$-relaxed strong edge coloring of $G$, thereby completing the proof.

Before starting the proof, we briefly describe the main idea.   
We suppose that there exists a $2$-optimal coloring  $\phi$ of $G$  containing bad edges.
Then, we characterize the structural properties of bad edges with respect to $\phi$ through a series of claims. Each claim is proved by contradiction: suppose the contrary, we  construct a new coloring that contradicts the $2$-optimality of $\phi$.
Finally, according to  the special structural properties of the bad edges, in the proof of \Cref{lemma:G notin Gp},
we construct a new coloring of $G$ based on $\phi$ by
recoloring certain edges of $G$. This new coloring is a good coloring with fewer bad edges, leading to  a contradiction to the $1$-optimality (and thus the $2$-optimality) of $\phi$, and so  any $2$-optimal coloring of $G$  contains no bad edges.


Let $\phi$ be  a $2$-optimal coloring of $G$. 
Suppose to the contrary that
 $\kappa_{1}(\phi)> 0$, that is, there exist bad edges (and hence bad pairs) with respect to $\phi$ in $G$. 
For brevity, we shall  refer to the abbreviation “the bad edges with respect to $\phi$” simply as  “bad edges”, and “the bad pairs with respect to $\phi$” as “bad pairs".

Note that every bad pair $e$ and $f$  in $G$ are $2$-neighbors of Type $6$ of each other.
We proceed to establish several structural  properties of  bad edges in $G$.

 

 \begin{claim}\label{claim:BE-type6-twice}
Let $e$ be a  bad edge of  $G$. For any color $\alpha \in [1,\Delta^{2}-1]\setminus \phi(F(e))$,  there are at least two edges in $T_{6}(e)$ being colored $\alpha$ in $\phi$. 
This implies that $|\phi(N^{2-}(e))|=\Delta^{2}-1$.
\end{claim}
\begin{proof}
Let $\alpha_{0}=\phi(e)$.
Since $e$ is a bad edge,
 $\alpha_{0}$ must appear on at least two edges in $T_{6}(e)$.
Suppose on the contrary that  there exists a color $\alpha \in [1,\Delta^{2}-1]\setminus \phi(F(e))$ such that  at most one edge in $T_{6}(e)$ is colored $\alpha$ under $\phi$.
We now  recolor the edge  $e$ with the color  $\alpha$.
This gives  a new coloring $\psi$ of $G$, which is obviously good.
It is easy to see that $\kappa_{1}(\psi,\alpha_{0})\le\kappa_{1}(\phi,\alpha_{0})-1$, $\kappa_{1}(\psi,\alpha)\le \kappa_{1}(\phi,\alpha)+1$, and $\kappa_{1}(\psi,\beta)=\kappa_{1}(\phi,\beta)$ for any color $\beta \in  [1,\Delta^{2}-1]\setminus\{\alpha_{0},\alpha\}$.
Hence, we have $\kappa_{1}(\psi)\le \kappa_{1}(\phi)$.
Moreover,  it holds that $\kappa_{2}(\psi,\alpha_{0})\le\kappa_{2}(\phi,\alpha_{0})-2$, $\kappa_{2}(\psi,\alpha)\le\kappa_{2}(\phi,\alpha)+1$, and $\kappa_{2}(\psi,\beta)=\kappa_{2}(\phi,\beta)$ for any color $\beta \in  [1,\Delta^{2}-1]\setminus\{\alpha_{0},\alpha\}$.
Therefore,  we  have $\kappa_{2}(\psi)<\kappa_{2}(\phi)$,  contradicting the $2$-optimality of $\phi$. The claim is proved.
\end{proof}

 \begin{claim}  \label{claim:BE-properties}\
 Let $e$ be a  bad edge of  $G$. Then, 
the following five properties hold:

$(1)$  $\phi(F(e))\cap \phi(T_{6}(e))=\emptyset$;

$(2)$  $|\phi(T_{6}(e))|= \frac{1}{2}|T_{6}(e)|$, (that is, the colors on the edges in $T_{6}(e)$ appear in pairs);

$(3)$ $|\phi(F(e))|=|F(e)|=\Delta^{2}-1-\frac{1}{2}|T_{6}(e)|$,   (that is,   all edges in $F(e)$ receive different colors);

$(4)$ $C_{e}^{\Delta}=\emptyset$, $($hence $T_{1}(e)=T_{2}(e)=T_{3}(e)=\emptyset$ and $N^{2}(e)=T_{4}(e)\cup T_{5}(e)\cup T_{6}(e)$$)$;

$(5)$ for any $w\in N(u)\cup N(v)$, $d(w)=\Delta$.
\end{claim}
\begin{proof}
By \Cref{claim:BE-type6-twice}, every color in $ [1,\Delta^{2}-1]\setminus \phi(F(e))$ appears on  at least two edges in $T_{6}(e)$.  This implies that $|\phi(T_{6}(e))\setminus\phi(F(e))|\le \frac{1}{2}|T_{6}(e)|$. 
Then, by
	 Observation \ref{Obser:shangjie}, we have
	 	\allowdisplaybreaks[4] 
		\begin{align*}	
	|\phi(N^{2-}(e))|	= &	|\phi(F(e))|+|\phi(T_{6}(e))\setminus\phi(F(e))|\\
		\le & \ 
	|F(e)|+ \frac{1}{2}|T_{6}(e)|\\
			\le & \
			\Delta^{2}-1-\frac{1}{2}|C_{e}^{\Delta}|-|T_{1}(e)|-\frac{1}{2}|T_{2}(e)|.
	\end{align*}
Recall that  $|\phi(N^{2-}(e))|=\Delta^{2}-1$ (see \Cref{claim:BE-type6-twice}),  
we must have
\allowdisplaybreaks[4] 
\begin{align}		|\phi(T_{6}(e))\setminus\phi(F(e))|	= & \	\frac{1}{2}|T_{6}(e)|,\\	
		|\phi(F(e))|= & \ 	|F(e)|,\\
		C_{e}^{\Delta} = & \ \emptyset, \\	
|F(e)|  = & \ \Delta^{2}-1-\frac{1}{2}|T_{6}(e)|.
\end{align}
Then the first four  properties are easy to see due to the above equations and \Cref{Obser:3-cycle}.
Finally, the last property
 follows directly from  Equation (6) and  \Cref{Obser:shangjie}.
\end{proof}


According to   Claim \ref{claim:BE-properties}(2), we immediately have the following observation.
\begin{observation}\label{Obser:bad edge-exactly2}
Every bad edge $e$ of $G$  has exactly two $2$-neighbors of Type $6$ with the same color as $e$.
\end{observation}

The next claim characterizes the properties of $1$-neighbors of a bad edge of $G$, which will be useful in the subsequent proofs.

 \begin{claim}\label{claim:BE-1-neighbor-type6-twice}
	Let $e$ be a  bad edge of  $G$. For any  $1$-neighbor $f$ of $e$, 
	if there exists some color $\alpha \in [1,\Delta^{2}-1]\setminus  \phi(F(f)\cup\{f\})$, then  there are at least two edges in $T_{6}(f)$ being colored $\alpha$ in $\phi$. 
	This implies that  $|\phi(N^{2-}[f])|=\Delta^{2}-1$ and  $|\phi(N^{2-}(f))|\ge \Delta^{2}-2$.
\end{claim}
\begin{proof}
	Let   $e$ be a bad edge  in $G$   and
 $f$ be a $1$-neighbor of $e$ with $[1,\Delta^{2}-1]\setminus  \phi(F(f)\cup\{f\})\neq \emptyset$. For convenience, let $\alpha_{0}=\phi(e)$ and $\alpha_{1}=\phi(f)$.	We prove by contradiction.
 Suppose that there is a color $\alpha\in [1,\Delta^{2}-1]\setminus  \phi(F(f)\cup\{f\})$ appears at most once on edges in $T_{6}(f)$. 
We recolor $f$ with the color $\alpha$ and $e$ with the color $\alpha_{1}$ to obtain a new  coloring $\psi$ of $G$. It is clear that $\alpha \notin \psi(F(f))$.	
	Since $f\in N(e)$, by	\Cref{claim:BE-properties}, we have $\alpha_{1} \notin \phi(N^{2-}(e)\setminus\{f\})$ and so $\alpha_{1} \notin \psi(N^{2-}(e))$.	
Hence,   $\psi$ is a good coloring of $G$.
It is easy to check that $\kappa_{1}(\psi,\alpha_{0})\le \kappa_{1}(\phi,\alpha_{0})-1$,  $\kappa_{1}(\psi,\alpha_{1})\le\kappa_{1}(\phi,\alpha_{1})$,  $\kappa_{1}(\psi,\alpha)\le\kappa_{1}(\phi,\alpha)+1$, and  $\kappa_{1}(\psi,\beta)=\kappa_{1}(\phi,\beta)$ for any color $\beta \in  [1,\Delta^{2}-1]\setminus\{\alpha_{0},\alpha_{1},\alpha\}$. 
It follows that $\kappa_{1}(\psi)\le\kappa_{1}(\phi)$.
Moreover,  since  $\alpha$ appears at most once on edges in $T_{6}(f)$,
we have
$\kappa_{2}(\psi,\alpha_{0})=\kappa_{2}(\phi,\alpha_{0})-2$,  $\kappa_{2}(\psi,\alpha_{1})\le\kappa_{2}(\phi,\alpha_{1})$,  $\kappa_{2}(\psi,\alpha)\le\kappa_{2}(\phi,\alpha)+1$, and  $\kappa_{2}(\psi,\beta)=\kappa_{2}(\phi,\beta)$ for any color $\beta \in  [1,\Delta^{2}-1]\setminus\{\alpha_{0},\alpha_{1},\alpha\}$. 
Therefore,  $\kappa_{2}(\psi)<\kappa_{2}(\phi)$. This contradicts the $2$-optimality of $\phi$.
\end{proof}

Recall that   every bad edge $e$ of $G$  has exactly two $2$-neighbors of Type $6$ with the same color as $e$ (see \Cref{Obser:bad edge-exactly2}),
we proceed to establish the structural property of  these two $2$-neighbors.

\begin{claim}  \label{claim:BE-P6}
Let  $e=uv$  be a bad edge of $G$, and  let 
 $e_{1}$ and $e_{2}$ be the unique two  $2$-neighbors of $e$ with $\phi(e_{1})=\phi(e_{2})=\phi(e)$.
Then, we have    $|\{e_{1},e_{2}\}\cap N_{u}^{2}(e)|=1$ and
$|\{e_{1},e_{2}\}\cap N_{v}^{2}(e)|=1$.
\end{claim}
\begin{proof}
	If not, by symmetry,
	 we may assume that $\{e_{1},e_{2}\}\subseteq N_{u}^{2}(e)$. 
	Denote by $f$ the edge in  $N(e)$ that is adjacent to to $e_{1}$. It is clear that  $e,e_{1}\in N(f)$,  $e_{2}\in N^{2}(f)$ and 	$e_{2}\notin C^{\Delta}_{f}$. Since $\phi(e)=\phi(e_{1})=\phi(e_{2})$, by \Cref{claim:BE-1-neighbor-type6-twice} and  \Cref{Obser:shangjie},  we  have 
		\allowdisplaybreaks[4] 
		\begin{align*}	
	|\phi(N^{2-}(f))|	\le & \ 	|F(f)|-1+\frac{1}{2}(|T_{6}(f)|-1)\\
	\le & \ 
	\Delta^{2}-\frac{5}{2}-\frac{1}{2}|C_{f}^{\Delta}|-|T_{1}(f)|-\frac{1}{2}|T_{2}(f)|.
\end{align*}
As $|\phi(N^{2-}(f))|$ is an  integer, $|\phi(N^{2-}(f))|\le 	\Delta^{2}-3$, contradicting the fact  that $|\phi(N^{2-}(f))|\ge\Delta^{2}-2$ (refer to  \Cref{claim:BE-1-neighbor-type6-twice}). This finishes the proof of the claim.	
\end{proof}

 \begin{claim}  \label{claim:BE-type5-empty}
 Let $e$ be a  bad edge $e$ of $G$. Then,
$T_{5}(e)=\emptyset$ and  $N^{2}(e)=T_{4}(e)\cup T_{6}(e)$.
	 This implies that, for any $f\in N(e)$, $C_{f}^{\Delta}=\emptyset$ $($and so $T_{1}(f)=T_{2}(f)=T_{3}(f)=\emptyset$ and $N^{2}(f)=  T_{4}(f)\cup T_{5}(f)\cup T_{6}(f)$$)$.
\end{claim}
\begin{proof}
	Suppose that   $T_{5}(e)\neq\emptyset$. 
	Without loss of generality, let $e=uv$ and  $g$ be an edge in $T_{5}(e)\cap N_{u}^{2}(e)$.
    Let  $f$ be an edge in $N_{u}(e)$ that is adjacent to $g$. 
    It follows that 
	 $|C_{f}^{\Delta}|\ge2$. 	
	By \Cref{claim:BE-P6}, there exists one edge $e_{1}$  in	$N_{u}^{2}(e)$ with $\phi(e_{1})=\phi(e)$.
	Notice that $e\in N(f)$ and $e_{1}\in N^{2-}(f)$,  according to   \Cref{claim:BE-1-neighbor-type6-twice}  and  \Cref{Obser:shangjie}, we immediately  have 
		\allowdisplaybreaks[4] 
		\begin{align*}	
	|\phi(N^{2-}(f))|	\le & \ 	|F(f)|+\frac{1}{2}(|T_{6}(f)|-1)\\
	\le & \ \Delta^{2}-1-\frac{1}{2}|C_{f}^{\Delta}|-|T_{1}(f)|-\frac{1}{2}|T_{2}(f)|-\frac{1}{2}\\
		\le & \ \Delta^{2}-1-\frac{1}{2}\times2-|T_{1}(f)|-\frac{1}{2}|T_{2}(f)|-\frac{1}{2}\\
	= & \
	\Delta^{2}-\frac{5}{2}-|T_{1}(f)|-\frac{1}{2}|T_{2}(f)|.
\end{align*}
	Since $|\phi(N^{2-}(f))|$ is an integer, we have $|\phi(N^{2-}(f))|\le 	\Delta^{2}-3$, contradicting  \Cref{claim:BE-1-neighbor-type6-twice}.
	Therefore,  we have $T_{5}(e)=\emptyset$. Then, by \Cref{claim:BE-properties}(4), we have $N^{2}(e)=T_{4}(e)\cup T_{6}(e)$.  The claim is proved.
\end{proof}

 \begin{claim}  \label{claim:BE-2-neighborhood}
 Let $e=uv$ be a bad edge of $G$. Then, 
  $|\phi(N(e)\cup N_{u}^{2}(e))|=|N(e)\cup N_{u}^{2}(e)|=\Delta^{2}-1$ and $|\phi(N(e)\cup N_{v}^{2}(e))|=|N(e)\cup N_{v}^{2}(e)|=\Delta^{2}-1$ $($refer to  \Cref{fig:BE}, all the bold edges receive different colors$)$.
	\begin{figure}[htbp]  
	\centering
	\resizebox{10cm}{2.5cm}{\includegraphics{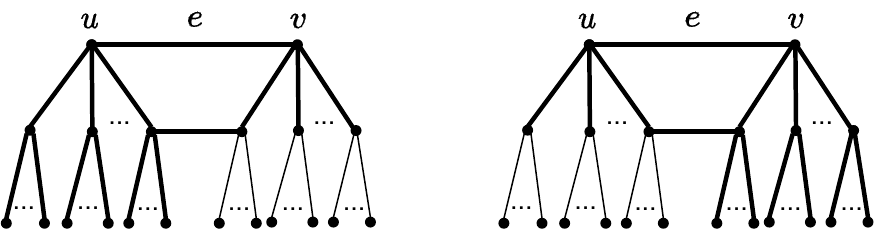}}
	\caption{The illustration of \Cref{claim:BE-2-neighborhood}.}
	\label{fig:BE}
\end{figure}
\end{claim}
\begin{proof}
Let $\alpha_{0}=\phi(e)$. Since any vertex in $N(u)\cup N(v)$ is of degree $\Delta$ (see \Cref{claim:BE-properties}(5)),
	   $|N(e)\cup N_{u}^{2}(e)|=|N(e)\cup N_{v}^{2}(e)|=\Delta^{2}-1$. In the following, we prove that $|\phi(N(e)\cup N_{u}^{2}(e))|= |\phi(N(e)\cup N_{v}^{2}(e))|=\Delta^{2}-1$.

       Suppose  that  $|\phi(N(e)\cup N_{u}^{2}(e))|< \Delta^{2}-1$. 	Let 
	   $\alpha \in  [1,\Delta^{2}-1]\setminus \phi(N(e)\cup N_{u}^{2}(e))$.   
	   According to \Cref{claim:BE-P6}, we may assume that
 $e_{1}\in N_{u}^{2}(e)\cap T_{6}(e)$ and $e_{2}\in N_{v}^{2}(e)\cap T_{6}(e)$ are  the  two distinct $2$-neighbors of $e$ being colored $\alpha_{0}$.   It is obvious that $\alpha\neq \alpha_{0}$.	
	 By \Cref{claim:BE-type5-empty}, $N^{2-}(e)=N(e)\cup T_{4}(e)\cup T_{6}(e)$.
	  Note that $T_{4}(e) = N_{u}^{2}(e)\cap N_{v}^{2}(e)$, we must have $\alpha \notin \phi(T_{4}(e))$ (as $\alpha \notin \phi(N_{u}^{2}(e))$).
     Because $|\phi(N^{2-}(e))|=\Delta^{2}-1$ (see \Cref{claim:BE-type6-twice}) and the colors on the edges in $T_{6}(e)$ appear in pairs (see Claim~\ref{claim:BE-properties}(2)),
      the color $\alpha$ appears on exactly two edges  in   $N_{v}^{2}(e)\cap T_{6}(e)$.
	  It follows that  there are exactly two edges $f_{1}$ and $f_{2}$ in $N_{u}^{2}(e)\cap T_{6}(e)$ such that $\phi(f_{1})=\phi(f_{2}) \neq\alpha$. 	 
	 Denote by $h_{1}$ and $h_{2}$ the edges in $N(e)$ being adjacent to $f_{1}$ and $f_{2}$, respectively. It is clear that $h_{1}\neq h_{2}$. 
We may assume that $e_{1}\notin N(h_{1})$ (and so $e_{1}\in N^{2}(h_{1})$). According to \Cref{claim:BE-type5-empty},   we have $C_{h_{1}}^{\Delta}=\emptyset$ and  $N^{2}(h_{1})=T_{4}(h_{1})\cup T_{5}(h_{1}) \cup T_{6}(h_{1})$.

We first  prove that  $e_{1}\in T_{6}(h_{1})$.
If  $e_{1}\in T_{5}(h_{1})$, then as $C_{h_{1}}^{\Delta}=\emptyset$, we must have $e_{1}\in T_{5}(e)$, contradicting the assumption  that $e_{1}\in  T_{6}(e)$.
If $e_{1}\in T_{4}(h_{1})$, note that 
$e,e_{1},f_{1}\in F(h_{1})$, $f_{2}\in N^{2}(h_{1})$, $\phi(e)=\phi(e_{1})$, and $\phi(f_{1})=\phi(f_{2})$, hence by \Cref{claim:BE-1-neighbor-type6-twice} and   \Cref{Obser:shangjie},
 we have
	\allowdisplaybreaks[4] 
		\begin{align*}	
	|\phi(N^{2-}(h_{1}))|	\le & \ 	|F(h_{1})|-1+\frac{1}{2}(|T_{6}(h_{1})|-1)\\
	\le & \
	\Delta^{2}-\frac{5}{2}-\frac{1}{2}|C_{h_{1}}^{\Delta}|-|T_{1}(h_{1})|-\frac{1}{2}|T_{2}(h_{1})|.
\end{align*}
	Since $|\phi(N^{2-}(h_{1}))|$ is an integer, $|\phi(N^{2-}(h_{1}))|\le 	\Delta^{2}-3$, contradicting the conclusion in  \Cref{claim:BE-1-neighbor-type6-twice} that  $|\phi(N^{2-}(h_{1}))|\ge \Delta^{2}-2$.
  Therefore,  we must have  $e_{1}\in T_{6}(h_{1})$.

    Next, we prove that there is no edge in $N^{2-}(h_{1})\setminus\{e,e_{1}\}$ being colored $\alpha_{0}$. 
     If not, let    $e^{*}$ be an edge in $N^{2-}(h_{1})\setminus\{e,e_{1}\}$ that is colored with $\alpha_{0}$. 
      Recall that  no edge  in $N^{2-}(e)$ is colored $\alpha_{0}$ except  $e_{1}$ and $e_{2}$.
    Since any edge in $(N(h_{1})\cup T_{4}(h_{1}))\setminus\{e\}$ is also in $N^{2-}_{u}(e)$, 
we must have $e^{*}\notin N(h_{1})\cup T_{4}(h_{1})$ and so  $e^{*} \in  T_{5}(h_{1}) \cup T_{6}(h_{1})$ (note that $N^{2-}(h_{1})=N(h_{1})\cup T_{4}(h_{1})\cup T_{5}(h_{1}) \cup T_{6}(h_{1})$). 
     If $e^{*} \in T_{5}(h_{1})$,
  notice that   $e,e^{*},f_{1}\in F(h_{1})$, $f_{2}\in N^{2}(h_{1})$, $e_{1}\in T_{6}(h_{1})$, $\phi(e)=\phi(e^{*})=\phi(e_{1})$ and $\phi(f_{1})=\phi(f_{2})$, hence  by  \Cref{claim:BE-1-neighbor-type6-twice} and   \Cref{Obser:shangjie},
   we have
   	\allowdisplaybreaks[4] 
  \begin{align*}	
  	|\phi(N^{2-}(h_{1}))|	\le & \ 	|F(h_{1})|-1+\frac{1}{2}(|T_{6}(h_{1})|-2)\\
  	\le & \
  	\Delta^{2}-3-\frac{1}{2}|C_{h_{1}}^{\Delta}|-|T_{1}(h_{1})|-\frac{1}{2}|T_{2}(h_{1})|,
  \end{align*}
this is again a contradiction to   \Cref{claim:BE-1-neighbor-type6-twice}. 
If $e^{*} \in T_{6}(h_{1})$, 
note that  $e,f_{1}\in N(h_{1})$, $f_{2}\in N^{2}(h_{1})$, $e_{1},e^{*}\in T_{6}(h_{1})$, $\phi(e)=\phi(e_{1})=\phi(e^{*})$ and $\phi(f_{1})=\phi(f_{2})$,  thus 
   	\allowdisplaybreaks[4] 
  		\begin{align*}	
  	|\phi(N^{2-}(h_{1}))|	\le & \ 	|F(h_{1})|+\frac{1}{2}(|T_{6}(h_{1})|-3)\\
  \le & \
  	\Delta^{2}-\frac{5}{2}-\frac{1}{2}|C_{h_{1}}^{\Delta}|-|T_{1}(h_{1})|-\frac{1}{2}|T_{2}(h_{1})|.
  \end{align*}
Again we have $|\phi(N^{2-}(h_{1}))|\le 	\Delta^{2}-3$,  a contradiction to   \Cref{claim:BE-1-neighbor-type6-twice}.

   Now we can exchange the colors of
  $e$ and $h_{1}$ in $\phi$ to get a new  coloring $\psi$ of $G$. 
  Since $h_{1}\in N(e)$, by \Cref{claim:BE-properties},  we have $\phi(h_{1})\notin\phi(N^{2-}(e)\setminus\{h_{1}\})$.
It follows that $\psi(e)=\phi(h_{1})\notin\psi(N^{2-}(e))$ and thus $e$ is not a bad edge with respect to $\psi$.
   Because  $e_{1}$ is the only  edge in $N^{2-}(h_{1})$ that is  colored with $\alpha_{0}$ in $\psi$ and $e_{1}\in T_{6}(h_{1})$,  hence we have $\psi(h_{1})=\alpha_0 \notin\psi(F(h_{1}))$ and so
   $h_{1}$ is not a bad edge with respect to $\psi$. 
  Therefore,  $\psi$ is a good coloring of $G$.
Moreover, it is easy to check that
  $\kappa_{1}(\psi,\alpha_{0})\le\kappa_{1}(\phi,\alpha_{0})-1$,  $\kappa_{1}(\psi,\phi(h_{1}))\le\kappa_{1}(\phi,\phi(h_{1}))$ and  $\kappa_{1}(\psi,\beta)=\kappa_{1}(\phi,\beta)$ for any color $\beta \in  [1,\Delta^{2}-1]\setminus\{\alpha_{0},\phi(h_{1})\}$. 
Hence,  $\kappa_{1}(\psi)<\kappa_{1}(\phi)$, which contradicts the $1$-optimality of $\phi$. 
  Consequently,  $|\phi(N(e)\cup N_{u}^{2}(e))|= \Delta^{2}-1$. By symmetry, we also have $|\phi(N(e)\cup N_{v}^{2}(e))|= \Delta^{2}-1$. 
  The claim is proved.
    \end{proof}

Based on Claims~\ref{claim:BE-properties}(2) and~\ref{claim:BE-2-neighborhood}, we immediately observe the following.

\begin{observation}\label{Obser:bad edge-type6}
Let $e=uv$ be a bad edge of $G$. Then,  $|T_{6}(e)|$ is even, and 
$$|\phi(T_{6}(e)\cap N^{2}_{u}(e))|=|T_{6}(e)\cap N^{2}_{u}(e)|= |T_{6}(e)\cap N^{2}_{v}(e)|=|\phi(T_{6}(e)\cap N^{2}_{v}(e))|.$$
\end{observation}

\begin{figure}[htbp]  
	\centering
	\resizebox{11.5cm}{2.7cm}{\includegraphics{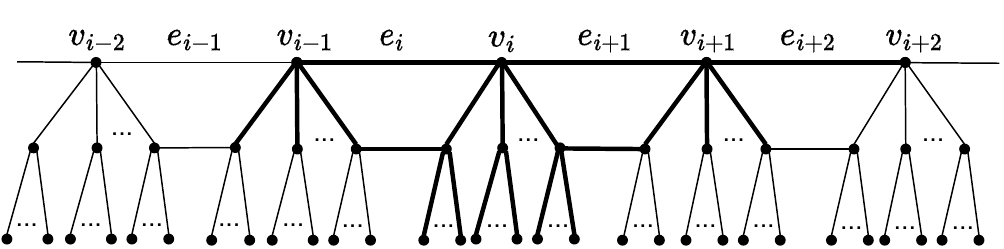}}
	\caption{The illustration of the  edge set $M_{e_{i}} \  (2\le i \le k)$.}
	\label{fig:Me}
\end{figure}

	Let $k\ge2$ be an integer. Suppose 
$P_{k}=v_{0}v_{1}v_{2}\ldots v_{k}$ is an  induced path in $G$,  with the edge $v_{0}v_{1}$ being a bad edge. In the following two claims, we  investigate   properties of this path.
For every $i\in [1,k]$,  let $e_{i}=v_{i-1}v_{i}$ and let
\begin{equation*}\label{eq:Me}
	\begin{array}{rl}
		M_{e_{i}}=\begin{cases}
			N(e_{1})\cup N_{v_{1}}^{2}(e_{1}), &  i=1,\\
			(N[e_{i}]\setminus \{e_{i-1}\})\cup N_{v_{i}}^{2}(e_{i}), & 2\le i \le k.
		\end{cases}
	\end{array}
\end{equation*}
As shown in \Cref{fig:Me}, the  edge set $M_{e_{i}}$ ($2\le i \le k$) is indicated by bold edges. Note that $e_{1}\notin M_{e_{1}}$, and that $e_{i}\in M_{e_{i}}$  for each $2\le i \le k$.

\begin{claim} \label{claim:path-1}
	Let $k\ge2$ be an integer. Suppose 
	$P_{k}=v_{0}v_{1}v_{2}\ldots v_{k}$ is an  induced path in $G$  with  $v_{0}v_{1}$ being a bad edge.
	Then,  for every  $1\le i\le k$,   $|\phi(M_{e_{i}})|=|M_{e_{i}}|=\Delta^{2}-1$.
\end{claim}
\begin{proof}
Because $e_{1}=v_{0}v_{1}$ is a bad edge, 	it follows from   \Cref{claim:BE-2-neighborhood} that this claim holds for $i=1$.
	 For any integer $ 2 \le i \le k$, it is clear  that	$|\phi(M_{e_{i}})|
	 \le |M_{e_{i}}|
	  = 	|N_{v_{i-1}}[e_{i}]\setminus \{e_{i-1}\}|+|N_{v_{i}}^{2-}(e_{i})|\le (\Delta-1)+ \Delta(\Delta-1)=
	\Delta^{2}-1$. Therefore, we just need to prove  $|\phi(M_{e_{i}})|=\Delta^{2}-1$ for every $ 2 \le i \le k$.
		For convenience, let  $\alpha_{i}=\phi(e_{i})$ for each $i\in[1,k]$.	
		We proceed by  induction on $k$.

	For $k=2$, if  $|\phi(M_{e_{2}})|<\Delta^{2}-1$, then we can recolor $e_{2}$ with some color $\alpha\in  [1,\Delta^{2}-1]\setminus \phi(M_{e_{2}})$ and $e_{1}$ with the color $\alpha_{2}$. This yields a new coloring of $G$ called $\psi$. 
	It is obvious  that $\alpha_{2}\neq\alpha$.
	
Since $e_{1}$ is a bad edge with respect to $\phi$ and $e_{2}\in N(e_{1})\subseteq F(e)$, by  Claims~\ref{claim:BE-properties} and~\ref{claim:BE-2-neighborhood}, we have $\alpha_{2}\notin \phi(N^{2-}(e_{1})\setminus\{{e_{2}}\})$ and so $\psi(e_{1})=\alpha_{2}\notin \psi(N^{2-}(e_{1}))$.  It follows that $\psi(e_{1})\notin \psi(F(e_{1}))$ and $e_1$ is not bad edge of $G$ with respect to $\psi$.
Recall that $|\phi(M_{e_{1}})|=|M_{e_{1}}|=\Delta^{2}-1$, there is exactly one edge $f$ in $M_{e_{1}}$ being colored $\alpha$ under $\phi$. 
Because  $\alpha\notin   \phi(M_{e_{2}})$, we must have $f\in M_{e_{1}}\setminus M_{e_{2}}$ and so $f$ is a $2$-neighbor of $e_{2}$. 
Note that $T_{4}(e_{2})\subseteq M_{e_{2}}$, hence $f\notin T_{4}(e_{2})$. 
Note also that $C_{e_{1}}^{\Delta}=\emptyset$ (see \Cref{claim:BE-properties}(4)) and $T_{5}(e_{1})=\emptyset$ (see \Cref{claim:BE-type5-empty}), we must have $f\in T_{6}(e_{2})$ and thus $\psi(e_{2})=\alpha \notin \psi(F(e_{2}))$. Therefore,   $\psi$ is a good coloring of $G$.

It is easy to see that
		$\kappa_{1}(\psi,\alpha_{1})\le\kappa_{1}(\phi,\alpha_{1})-1$, 	$\kappa_{1}(\psi,\alpha_{2})\le\kappa_{1}(\phi,\alpha_{2})$,
	$\kappa_{1}(\psi,\alpha)\le\kappa_{1}(\phi,\alpha)+1$, and  $\kappa_{1}(\psi,\beta)=\kappa_{1}(\phi,\beta)$ for any color $\beta \in  [1,\Delta^{2}-1]\setminus\{\alpha_{1},\alpha_{2},\alpha\}$. It follows that $\kappa_{1}(\psi)\le\kappa_{1}(\phi)$.
Moreover,   it is clear that 
	$\kappa_{2}(\psi,\alpha_{1})=\kappa_{2}(\phi,\alpha_{1})-2$,  $\kappa_{2}(\psi,\alpha_{2})\le\kappa_{2}(\phi,\alpha_{2})$,  $\kappa_{2}(\psi,\alpha)\le\kappa_{2}(\phi,\alpha)+1$, and  $\kappa_{2}(\psi,\beta)=\kappa_{2}(\phi,\beta)$ for any color $\beta \in  [1,\Delta^{2}-1]\setminus\{\alpha_{1},\alpha_{2},\alpha\}$. 
This implies that $\kappa_{2}(\psi)<\kappa_{2}(\phi)$, which is a  contradiction to the $2$-optimality of $\phi$. Hence, we have
	 $|\phi(M_{e_{2}})|=|M_{e_{2}}|=\Delta^{2}-1$.

Next, we  consider the case  $k\ge 3$.
	Assume that   $|\phi(M_{e_{i}})|=|M_{e_{i}}|=\Delta^{2}-1$ holds for any integer $1\le i\le k-1$. In what follows, we prove $|\phi(M_{e_{k}})|=\Delta^{2}-1$.
If not,  let $\alpha\in  [1,\Delta^{2}-1]\setminus \phi(M_{e_{k}})$.	
It is obvious that $\alpha\neq \alpha_{k}$.  (Note that it is possible that  $\alpha= \alpha_{k-1}$.)
Now, we  recolor $e_{i}$ with the color $\alpha_{i+1}$ for each $i\in[1,k-1]$ and  $e_{k}$ with the color $\alpha$. This results in  a new coloring of $G$ called $\psi$, in which $\psi(e_{k})=\alpha$ and $\psi(e_{i})=\alpha_{i+1}$ for each $i\in[1,k-1]$.

Because   $P_{k}=v_{0}v_{1}v_{2}\ldots v_{k}$ is an  induced path in $G$  with $e_{1}$ being a bad edge    and  $|\phi(M_{e_{i}})|=|M_{e_{i}}|=\Delta^{2}-1$ for every $1\le i \le k-1$, we immediately observe the following.
\begin{observation}\label{ob:A}
 For any $1\le i\le k-2$, there are $e_{i+2}\in T_{6}(e_{i})$ and  $e_{i}\in T_{6}(e_{i+2})$. Moreover, if $k\ge 4$, then $\alpha_{i}\neq \alpha_{i+2}$ for any $2\le i\le k-2$.
\end{observation}

Recall that $|M_{e_{i}}|=\Delta^{2}-1$ for every $1\le i \le k-1$, the following observation follows directly.
\begin{observation}\label{ob:B}
	For any $1\le i \le k-1$,  we have $C_{e_{i}}^{\Delta}=\emptyset$ and  $T_{5}(e_{i})=\emptyset$. 
\end{observation}

Before proceeding with the proof,
we make two other useful observations.

\begin{observation}\label{ob:C}
The edge	$e_{k}$ has exactly one $2^-$-neighbor $f$ being colored  $\alpha$ under $\psi$. Moreover, $f\in (M_{e_{k-1}}\setminus M_{e_{k}})\cup\{e_{k-2}\}$  and $f\in T_{6}(e_{k})$. (This implies that, $e_{k}$ is not a bad edge with respect to $\psi$; however, it is possible that $f$ is  a  bad edge with respect to $\psi$.)
\end{observation}
\begin{proof}
Since $|\phi(M_{e_{k-1}})|=|M_{e_{k-1}}|=\Delta^{2}-1$ and  $\alpha\notin   \phi(M_{e_{k}})$,
	there is exactly one   edge  in $M_{e_{k-1}}\setminus M_{e_{k}}$ being colored $\alpha$ under $\phi$.
    Note that  $N^{2-}(e_{k})=\{e_{k-2}\}\cup M_{e_{k-1}}\cup M_{e_{k}}$, and that $\psi(e_{k-2})=\alpha_{k-1}$ (possibly $\alpha=\alpha_{k-1}$).
Therefore,  exactly one edge $f$ in $(M_{e_{k-1}}\setminus M_{e_{k}})\cup\{e_{k-2}\}$ (hence, exactly one  $2^-$-neighbor $f$ of $e_k$)  is colored $\alpha$ under $\psi$.
	Now, if  $f=e_{k-2}$,  then  it is clear that $f\in T_{6}(e_{k})$ due to  Observation~\ref{ob:A}. 
	If $f\in M_{e_{k-1}}\setminus M_{e_{k}}$,  then since $C_{e_{k-1}}^{\Delta}=\emptyset$ and $T_{5}(e_{k-1})=\emptyset$ (see Observation~\ref{ob:B}), it is easy to check that $f\in T_{6}(e_{k})$.
\end{proof}

\begin{observation}\label{ob:D}
For every $2\le i\le k$ with $i\neq k-1$, it holds that $\alpha_{i}=\psi(e_{i-1})\notin \psi(N^{2-}(e_{i-1}))$. In particular,  we have	$\alpha_{k-1}=\psi(e_{k-2})\notin \psi(N^{2-}(e_{k-2})\setminus\{e_{k}\})$.
(This, together with $e_k\in T_6(e_{k-2})$ (see Observation~\ref{ob:A}),  implies that every edge $e_i$ ($1\le i\le k-1$) is not a bad edge with respect to $\psi$.)
\end{observation}
\begin{proof}
	 We first prove that $\alpha_{k}\notin \psi(N^{2-}(e_{k-1}))$.
     Note that $e_k$ belongs to both $M_{e_{k-2}}$ and $M_{e_{k-1}}$. 
	 Since    
     $|\phi(M_{e_{k-2}})|=|M_{e_{k-2}}|=\Delta^{2}-1$ and  $|\phi(M_{e_{k-1}})|=|M_{e_{k-1}}|=\Delta^{2}-1$, 
     we must have $\alpha_{k}\notin \phi(M_{e_{k-2}}\setminus\{e_{k}\})$ and 
	 $\alpha_{k}\notin \phi(M_{e_{k-1}}\setminus\{e_{k}\})$.
     Recall that $\psi(e_{k-2})=\alpha_{k-1}$, $\psi(e_{k-1})=\alpha_{k}$, $\alpha_k\neq \alpha_{k-1}$, and  $\alpha_k\neq \alpha$.
	It follows that $\alpha_{k}\notin \psi(N^{2-}(e_{k-1}))$.
	
	Then, we prove that  for every $2\le i \le k-1$,	it holds that $\alpha_{i}\notin \psi(N^{2-}(e_{i-1})\setminus\{e_{i+1}\})$.
	 Recall that	 $e_{1}$ is a bad edge with respect to $\phi$ and $e_{2}\in N(e_{1})$, hence by  \Cref{claim:BE-2-neighborhood}, we have $\alpha_{2}\notin \phi(N^{2-}(e_{1})\setminus\{{e_{2}}\})$.
     It follows that $\alpha_{2}\notin \psi(N^{2-}(e_{1})\setminus\{e_{3}\})$ (note that possibly $\alpha_{2}=\alpha$).
	 Now if $k=3$, the proof is complete.
     Assume that $k\ge 4$.
 For every $3\le i \le k-1$, 	 because  
	 $|\phi(M_{e_{i-2}})|=|M_{e_{i-2}}|=\Delta^{2}-1$ and  $|\phi(M_{e_{i-1}})|=|M_{e_{i-1}}|=\Delta^{2}-1$, 
     we  have 
     $\alpha_{i}=\phi(e_{i})\notin \phi(M_{e_{i-2}}\setminus\{e_{i}\})$ and 
	 $\alpha_{i}=\phi(e_{i})\notin \phi(M_{e_{i-1}}\setminus\{e_{i}\})$.
	  It follows that   $\alpha_{i}\notin \psi(N^{2-}(e_{i-1})\setminus\{e_{i+1}\})$ for every $2\le i \le k-1$. 
	 
	Finally,  
	 due to Observation \ref{ob:A}, for every $2\le i \le k-2$, it holds that
	 $\alpha_{i}\neq \alpha_{i+2}= \psi(e_{i+1})$. This, together with  $\alpha_{i}\notin \psi(N^{2-}(e_{i-1})\setminus\{e_{i+1}\})$, implies that
	$\alpha_{i}\notin \psi(N^{2-}(e_{i-1}))$ for any $2\le i \le k-2$. 
	 Therefore, the observation is proved.
\end{proof}

In light of  Observations~\ref{ob:C} and~\ref{ob:D}, it follows easily that $\psi(e_{i}) \notin \psi(F(e_{i}))$ for every $1\le i \le k$. 
Thus, $\psi$ is a good coloring of $G$. 
Moreover,  these two observations  also imply that   $\kappa_{1}(\psi,\alpha_{i})\le\kappa_{1}(\phi,\alpha_{i})$ for every $2\le i \le k$ and  $\kappa_{1}(\psi,\alpha)\le\kappa_{1}(\phi,\alpha)+1$. 
Notice that  $\kappa_{1}(\psi,\alpha_{1})\le\kappa_{1}(\phi,\alpha_{1})-1$  and   $\kappa_{1}(\psi,\beta)=\kappa_{1}(\phi,\beta)$ for any color $\beta \in  [1,\Delta^{2}-1]\setminus\{\alpha_{1},\alpha_{2},\ldots,\alpha_{k},\alpha\}$,  we have $\kappa_{1}(\psi)\le\kappa_{1}(\phi)$.
It is easy to check that $\kappa_{2}(\psi,\alpha_{1})\le\kappa_{2}(\phi,\alpha_{1})-2$, $\kappa_{2}(\psi,\alpha)\le\kappa_{2}(\phi,\alpha)+1$,
and
$\kappa_{2}(\psi,\beta)\le\kappa_{2}(\phi,\beta)$ for any color $\beta \in  [1,\Delta^{2}-1]\setminus\{\alpha_{1},\alpha\}$. Therefore, we obtain $\kappa_{2}(\psi)<\kappa_{2}(\phi)$, a contradiction to the $2$-optimality of $\phi$. 
Hence, it holds that
 $|\phi(M_{e_{k}})|=|M_{e_{k}}|=\Delta^{2}-1$.
The claim is proved.
\end{proof}

\begin{claim} \label{claim:path-2}
Let   $k\ge2$ be an integer. Suppose 
$P_{k}=v_{0}v_{1}v_{2}\ldots v_{k}$ is an  induced path in $G$,  where $v_{0}v_{1}$ is a bad edge and $e_i=v_{i-1}v_i$ for each $i\in[1,k]$. 
Then  we have the following three conclusions.

	$(1)$ For every  $1\le i\le k$,  it holds that $C_{e_{i}}^{\Delta}=\emptyset$,  $T_{5}(e_{i})=\emptyset$  and  $N^{2}(e_{i})=T_{4}(e_{i})\cup T_{6}(e_{i})$;

	$(2)$ 
	For every $2\le i\le k$ and  $i\neq3$,  there is $\phi(e_{i})\notin \phi(N^{2-}(e_{i}))$; and while if $k\ge3$, then $\phi(e_{3})\notin \phi(N^{2-}(e_{3})\setminus\{e_{1}\})$;

	$(3)$ If $k\ge 3$, then  for every  $3\le i\le k$, there exists exactly one  edge  $h_{i}\in  N^{2}_{v_{i}}(e_{i})$ such that $\phi(h_{i})=\phi(e_{i-1})$; moreover,  it holds that $h_{i}\in T_{6}(e_{i})$.
\end{claim}
\begin{proof}	
	The first  conclusion holds for $i=1$ due to \Cref{claim:BE-type5-empty}.
	As $P_{k}=v_{0}v_{1}v_{2}\ldots v_{k}$ is an  induced path   with $e_{1}$ being a bad edge,  it follows from \Cref{claim:path-1} that
$|\phi(M_{e_{i}})|=|M_{e_{i}}|=\Delta^{2}-1$  for  any $1\le i\le k$. 
	Thus, for every  $2\le i\le k$, we must have  $C_{e_{i}}^{\Delta}=\emptyset$ and $T_{5}(e_{i})=\emptyset$ as otherwise there is a contradiction to the fact that $|M_{e_{i-1}}|=|M_{e_{i}}|=\Delta^{2}-1$.
Hence, conclusion  (1) is correct.

		Then, we prove  conclusion (2).
	Because
$|\phi(M_{e_{1}})|=|M_{e_{1}}|=\Delta^{2}-1$, $|\phi(M_{e_{2}})|=|M_{e_{2}}|=\Delta^{2}-1$ and $e_{2}\in M_{e_{1}}\cap M_{e_{2}}$, 
we must have	$\phi(e_{2})\notin \phi(N^{2-}(e_{2})\setminus\{e_{1}\})$.
	Notice that $e_{2}\in N(e_{1})$, and thus  $\phi(e_{2})\neq\phi(e_{1})$. Hence, 	$\phi(e_{2})\notin \phi(N^{2-}(e_{2}))$.
	While if $k\ge3$, 
	for every $3\le i \le k$, since
		$|\phi(M_{e_{i-1}})|=|M_{e_{i-1}}|=\Delta^{2}-1$,  $|\phi(M_{e_{i}})|=|M_{e_{i}}|=\Delta^{2}-1$ and $e_{i}\in M_{e_{i-1}}\cap M_{e_{i}}$,
	it holds that
		$\phi(e_{i})\notin \phi(N^{2-}(e_{i})\setminus\{e_{i-2}\})$.
	When $k\ge4$, for every $4\le i \le k$, since $e_{i-2},e_{i}\in M_{e_{i-2}}$ and $|\phi(M_{e_{i-2}})|=|M_{e_{i-2}}|=\Delta^{2}-1$, we have $\phi(e_{i})\neq\phi(e_{i-2})$ and so 	$\phi(e_{i})\notin \phi(N^{2-}(e_{i}))$.

Finally, we prove that  conclusion (3) is also true. 
  For every $3\le i\le k$,   since $|\phi(M_{e_{i}})|=|M_{e_{i}}|=\Delta^{2}-1$ and $e_{i-1}\notin M_{e_{i}}$,  there is  exactly one edge $h_{i}$ in $M_{e_{i}}$ such that $\phi(h_{i})=\phi(e_{i-1})$.
  Because $|\phi(M_{e_{i-1}})|=|M_{e_{i-1}}|=\Delta^{2}-1$ and $e_{i-1}\in M_{e_{i-1}}$, we have $\phi(e_{i-1})\notin \phi(M_{e_{i-1}}\setminus\{e_{i-1}\})$.
  Recall that $M_{e_{i}}=(N[e_{i}]\setminus \{e_{i-1}\})\cup N^{2}_{v_{i}}(e_{i})$,
  we must have  $h_{i}\in  N^{2}_{v_{i}}(e_{i})$
  as $N[e_{i}]\setminus \{e_{i-1}\}\subseteq M_{e_{i-1}}\setminus\{e_{i-1}\}$.
Due to  conclusion (1), there is 
 $N^{2}(e_{i})=T_{4}(e_{i})\cup T_{6}(e_{i})$. If $h_{i}\in T_{4}(e_{i})$, then $h_{i}\in M_{e_{i-1}}\setminus\{e_{i-1}\}$, which is a contradiction since $\phi(e_{i-1})\notin  \phi(M_{e_{i-1}}\setminus\{e_{i-1}\})$.
Therefore, we must have   $h_{i}\in T_{6}(e_{i})$ for every $3\le i\le k$.
This finishes the proof.
\end{proof}

\begin{claim}\label{claim:type6}
Let   $e=uv$ be a  bad edge  of $G$. 
  Then,  there are two vertices $u'\in N(u)$ and $v'\in N(v)$ such that  $u'v'\notin E(G)$, and $|N(uu')\cap T_{6}(e)|\ge2$ or $|N(vv')\cap T_{6}(e)|\ge2$. 
\end{claim}
\begin{proof}
Let $\alpha_0=\phi(e)$.
 According to \Cref{claim:BE-P6}, we may assume that
$e_{1}=u_{1}x_{1}\in N_{u}^{2}(e)\cap T_{6}(e)$ and $e_{2}=v_{1}y_{1}\in N_{v}^{2}(e)\cap T_{6}(e)$ are the  two distinct $2$-neighbors of $e$ being colored $\alpha_{0}$, where $u_{1}\in N(u)$ and  $v_{1}\in N(v)$.
	Denote by $f_{1}$  and $f_{2}$ the two edges $uu_{1}$ and  $vv_{1}$, respectively.
	For brevity,  let $\alpha_{1}=\phi(f_{1})$ and  $\alpha_{2}=\phi(f_{2})$. 
	By \Cref{claim:BE-properties}(3), it is clear that $\alpha_{1}\neq \alpha_{2}$.
	
	First we prove that   $T_{6}(e)\setminus \{e_{1},e_{2}\}\neq \emptyset$.
	Suppose on the  contrary that  $T_{6}(e)\setminus \{e_{1},e_{2}\}=\emptyset$.  
    Recall that $N^{2}(e)=T_{4}(e)\cup T_{6}(e)$ (see \Cref{claim:BE-type5-empty}), and thus $F(e)=N(e)\cup T_{4}(e)$ and $N^{2-}(e)=F(e)\cup \{e_{1},e_{2}\}$.
    Since $T_{6}(e)=\{e_{1},e_{2}\}$, by \Cref{claim:BE-properties}(3) and (5),  we have 
	$|F(e)|=|N(e)|+|T_{4}(e)|=\Delta^{2}-1-\frac{1}{2}|T_{6}(e)|=(\Delta^{2}-1)-1$ and $|N(e)|=2(\Delta-1)$.
	 Therefore,
$|T_{4}(e)|=(\Delta-1)^{2}-1$. 
This implies that,
for every $u'\in N(u)\setminus\{v\}$ and every
$v'\in N(v)\setminus\{u\}$, we must have $u'v'\in E(G)$ except when $u'=u_{1}$ and $v'=v_{1}$.
In other words, $\{e_{1},e_{2}\}$ is an edge cut of $G$. Refer to \Cref{fig:type6}.		
		\begin{figure}[htbp]  
			\centering
			\resizebox{15cm}{6cm}{\includegraphics{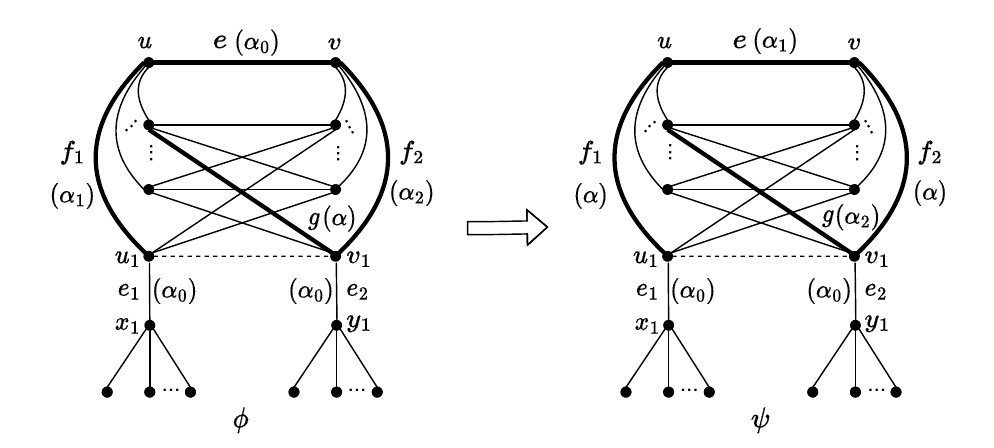}}
			\caption{The illustration of 
            the case   $T_{6}(e)=\{e_{1},e_{2}\}$ in \Cref{claim:type6}.}
			\label{fig:type6}
		\end{figure}

	Since $\phi(e_{1})=\phi(e_{2})=\alpha_0$ and $\Delta\ge3$, we have  
      	\allowdisplaybreaks[4] 
  		\begin{align*}	
 |\phi(N_{x_{1}}[e_{1}]\cup N_{y_{1}}[e_{2}]\cup \{f_{1},f_{2}\})|	
  \le & \ |N_{x_{1}}[e_{1}]|+|N_{y_{1}}[e_{2}]|-1+2\\
  	  \le & \  \Delta+\Delta-1+2\\
        <& \ \Delta^{2}-1.
  \end{align*}  
 It follows that there exists some color $\alpha$ in $[1,\Delta^{2}-1]\setminus\phi(N_{x_{1}}[e_{1}]\cup N_{y_{1}}[e_{2}]\cup \{f_{1},f_{2}\})$. (Note that $\alpha\notin\{\alpha_0,\alpha_1.\alpha_2\}$).
Recall that	$|\phi(N^{2-}(e))|=\Delta^{2}-1$ (see \Cref{claim:BE-type6-twice}), $N^{2-}(e)=F(e)\cup \{e_{1},e_{2}\}$, and $\phi(e_1)=\phi(e_2)=\alpha_0\neq\alpha$.
This, together with  \Cref{claim:BE-properties}(3), implies that there is exactly one edge $g$ in $F(e)=N(e)\cup T_{4}(e)$ being colored $\alpha$. 
Note that $|N(g)\cap\{e_{1},e_{2}\}|\le 1$, we may assume that $e_{1}\notin N(g)$.	
	Now, a new coloring $\psi$ can be obtained by recoloring $f_{1}$ and $f_{2}$ with the same color $\alpha$, $g$ with $\alpha_{2}$ and   $e$ with $\alpha_{1}$. 
	It is easy to check that $\psi$ is a good coloring of $G$.
	Moreover,  
it is straightforward to check that 
$\kappa_{1}(\psi,\alpha_{0})=\kappa_{1}(\phi,\alpha_{0})-1$ and  $\kappa_{1}(\psi,\beta)
\le\kappa_{1}(\phi,\beta)$ for any color $\beta \in  [1,\Delta^{2}-1]\setminus\{\alpha_{0}\}$. 
Therefore, we have  $\kappa_{1}(\psi)<\kappa_{1}(\phi)$,
contradicting the $1$-optimality of $\phi$. Thus  we must have $T_{6}(e)\setminus \{e_{1},e_{2}\}\neq \emptyset$.

It follows from $T_{6}(e)\setminus \{e_{1},e_{2}\}\neq \emptyset$  and  $|T_{6}(e)\cap N^{2}_{u}(e)|= |T_{6}(e)\cap N^{2}_{v}(e)|$ (see Observation \ref{Obser:bad edge-type6}) that there exists one vertex $u'\in N(u)$
such that $(N(uu')\cap T_{6}(e))\setminus\{e_{1}\}\neq \emptyset$. Without loss of generality, let $e_{3}=u'x'\in (N(uu')\cap T_{6}(e))\setminus\{e_{1}\}$.
Since $d(u')=\Delta$ (see  \Cref{claim:BE-properties}(5)) and $e_{3}=u'x'\in T_{6}(e)$,  we must have $|N(u')\cap (N(v)\setminus\{u\})|\le \Delta-2$. Therefore, there exists one vertex $v'\in N(v)$ such that $u'v'\notin E(G)$.
Hence, we also have $|N(v')\cap (N(u)\setminus\{v\})|\le \Delta-2$.
This implies that there is one edge $e_{4}=v'y'\in N(vv')\cap T_{6}(e)$.
We may assume that $N(uu')\cap T_{6}(e)=\{e_{3}\}$ and $N(vv')\cap T_{6}(e)=\{e_{4}\}$ as otherwise we are done.
It follows that $u'\neq u_{1}$,  $u'p\in E(G)$ for any $p\in N(v)\setminus\{v'\}$ and $v'q\in E(G)$ for any $q\in N(u)\setminus\{u'\}$.

A similar argument as above shows that there exists one vertex $v''\in N(v)$ such that $u_{1}v''\notin E(G)$ and there is one edge $e''=v''y''\in N(vv'')\cap T_{6}(e)$.
Notice that $u_1v'\in E(G)$ (as $u'\neq u_{1}$ and $v'q\in E(G)$ for any $q\in N(u)\setminus\{u'\}$), 
and thus $v''\neq v'$.
 We may also assume that
 $N(uu_{1})\cap T_{6}(e)=\{e_{1}\}$ and $N(vv'')\cap T_{6}(e)=\{e''\}$  as otherwise we are done.  
This implies that,  $u_{1}p\in E(G)$ for any $p\in N(v)\setminus\{v''\}$ and $v''q\in E(G)$ for any $q\in N(u)\setminus\{u_{1}\}$. 
 It is possible that $e_{2}\in \{e_{4},e''\}$. However, this will not affect the following arguments.

 For convenience, let $\alpha_{3}=\phi(uu')$, 
$\alpha_{4}=\phi(vv')$ and $\alpha_{5}=\phi(vv'')$.    
Since  $Q_{1}=vuu'$ is an induced path with $vu$ being a bad edge in $G$ and $uu_{1},vv''\in N(uu')\cup N^{2}_{u'}(uu')$,    by  \Cref{claim:path-1}, we have $\alpha_{1},\alpha_{5}\notin \phi (N_{x'}[e_{3}])$.
	Similarly, it holds that $\alpha_{1},\alpha_{5}\notin\phi ( N_{y'}[e_{4}])$ as $Q_{2}=uvv'$   is an  induced path in $G$.
Since  both $Q_{3}=vuu_{1}$ and $Q_{4}=uvv''$ are  induced paths in $G$,  we have $\alpha_{3},\alpha_{4}\notin \phi (N_{x_{1}}[e_{1}]) \cup \phi ( N_{y''}[e''])$. 
Now, we  recolor $e$ with $\alpha_{1}$, $uu'$ and $vv'$ with the same color $\alpha_{5}$,  $uu_{1}$  with  $\alpha_{3}$ and 
$vv''$ with  $\alpha_{4}$. This  gives rise to a new coloring called $\sigma$.
Note that every color in $\{\alpha_1,\alpha_3,\alpha_4,\alpha_5\}$ occurs exactly once in $N^{2-}(e)$ under $\phi$.
Recall that $u'v'\notin E(G)$, hence the two edges $uu'$ and $vv'$ are $2$-neighbors of Type 6 of each other.
It is not difficult to check that $\sigma$ is a good coloring of $G$ and $\kappa_{1}(\sigma)< \kappa_{1}(\phi)$, a contradiction to the $1$-optimality of $\phi$.

Therefore,  there exist two vertices $u'\in N(v)$ and $v'\in N(v)$ such that  $u'v'\notin E(G)$ and $|N(uu')\cap T_{6}(e)|\ge2$ or $|N(vv')\cap T_{6}(e)|\ge2$. 
This claim is proved.
\end{proof}

We have now all the ingredients to prove the case of    $G\notin  \mathcal{G}_{\Delta}$.

\begin{lemma}\label{lemma:G notin Gp}
	Let $G$ be a graph with maximum degree $\Delta\ge3$   that is not isomorphic to $K_{\Delta,\Delta}$.  	If $G\notin  \mathcal{G}_{\Delta}$,  then $\chi'_{ss}(G)\le \Delta^{2}-1$  and $\chi'_{(0,1)}(G)\le \Delta^{2}-1$. 
\end{lemma}
\begin{proof}
  Let $\phi$ be  a $2$-optimal coloring of $G$. By Observation \ref{Obser:no bad edge}, it  suffices to show that $G$ has no bad edge with respect to $\phi$.
  Suppose to the contrary, let  $e_{0}=uv$ be a bad edge with respect to $\phi$ in  $G$. Let $\alpha_{0}=\phi(e_{0})$.
   By \Cref{claim:BE-P6}, let $e_{1}\in N_{u}^{2}(e_{0})\cap T_{6}(e_{0})$ and $e_{2}\in N_{v}^{2}(e_{0})\cap T_{6}(e_{0})$ be  the unique two $2$-neighbors of $e_{0}$ being colored $\alpha_{0}$.
   
  According to \Cref{claim:type6}, there are  two vertices $u_{1}\in N(u)$ and $v_{1}\in N(v)$ such that $u_{1}v_{1}\notin E(G)$ and $|N(uu_{1})\cap T_{6}(e_{0})|\ge2$ or $|N(vv_{1})\cap T_{6}(e_{0})|\ge2$. 
   Without loss of generality, we may assume that $|N(vv_{1})\cap T_{6}(e_{0})|\ge2$.
    Denote by $g_{1}$ and $g_{2}$ the two edges $uu_{1}$ and $vv_{1}$, respectively.  Since $u_{1}v_{1}\notin E(G)$,   $g_{1}$ and $g_{2}$ do not lie on a common $4$-cycle in $G$.
    
  For brevity,  let  
   $\beta_{1}=\phi(g_{1})$ and $\beta_{2}=\phi(g_{2})$. 
   Since $C_{e_{0}}^{\Delta}=\emptyset$ (refer to \Cref{claim:BE-properties}(4)),   $Q_{1}=vuu_{1}$ and $Q_{2}=uvv_{1}$ are two induced paths in $G$ with $e_{0}$ being a bad edge.
     By  applying \Cref{claim:path-1} to $Q_{1}=vuu_{1}$ (resp.  $Q_{2}=uvv_{1}$),  there is exactly  one edge $h_{1}=s_{1}t_{1}$  in $N^{2}_{u_{1}}(g_{1})$   with $\phi(h_{1})=\beta_{2}$ (resp.  $h_{2}=s_{2}t_{2}$  in $N^{2}_{v_{1}}(g_{2})$  with $\phi(h_{2})=\beta_{1}$).
     Assume that  $s_{1}\in N(u_{1})$ and  $s_{2}\in N(v_{1})$.
    Refer to \Cref{fig:phi} for the illustration of the coloring $\phi$.
    We proceed by proving the following  claim.
   \begin{claim} \label{claim:h1h2}   	
   	$(1)$ $h_{1}\notin N^{2-}(e_{0})$, $h_{1}\notin N^{2-}(g_{2})$  and  	$h_{1}\in T_{6}(g_{1})$;
   	
 \qquad \ \  \ \	$(2)$ $h_{2}\notin N^{2-}(e_{0})$,   $h_{2}\notin N^{2-}(g_{1})$ and 
   	$h_{2}\in T_{6}(g_{2})$.
   \end{claim}
   \begin{proof}
   	By symmetry,  we only prove  (1) here.
   	Because $u_{1}v_{1}\notin E(G)$, we have $h_{1}\neq g_{2}$. It follows from \Cref{claim:BE-properties} that 
   	$h_{1}\notin N^{2-}(e_{0})$ since $\phi(h_{1})=\phi(g_{2})=\beta_{2}$ and $g_{2}\in N(e_{0})$.
   	Recall that $Q_{2}=uvv_{1}$ is an  induced path in $G$ with  $uv$ being a bad edge,
   	by  \Cref{claim:path-2}(2), 
   	we have  $\phi(g_{2})=\beta_{2}\notin \phi(N^{2-}(g_{2}))$.
  This, together with $\phi(h_{1})=\beta_{2}$, implies that $h_{1}\notin N^{2-}(g_{2})$.
  As $Q_{1}=vuu_{1}$ is an induced path with $vu$ being a bad edge,  by  \Cref{claim:path-2}(1), we have $N^{2}(g_{1})=T_{4}(g_{1})\cup T_{6}(g_{1})$. 
   	Since $h_{1}\in N^{2}_{u_{1}}(g_{1})$,  
   	$h_{1}$ is either in $T_{4}(g_{1})$ or  in $T_{6}(g_{1})$. 
   	If $h_{1}\in T_{4}(g_{1})$, then $h_{1}\in N^{2}(e_{0})$, which is a contradiction to the fact that 	 $h_{1}\notin N^{2-}(e_{0})$. 
   	Therefore, we must have  $h_{1}\in T_{6}(g_{1})$. This claim is true.
   \end{proof}
    	\begin{figure}[htbp]  
	\centering
	\resizebox{8cm}{5.6cm}{\includegraphics{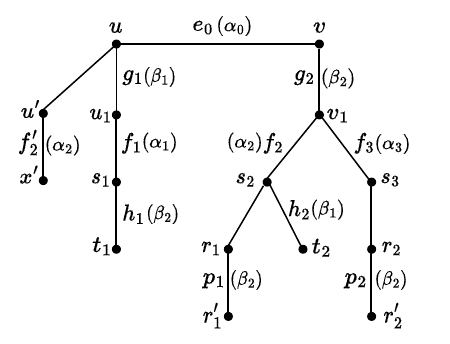}}
	\caption{The illustration of the coloring $\phi$}
	\label{fig:phi}
\end{figure} 
   
 We use  $f_{1}$ and  $f_{2}$ to denote  the two edges  $u_{1}s_{1}$ and  $v_{1}s_{2}$, respectively. 
   It is obvious that $f_{1},f_{2}\in N^{2}(e_{0})$. 
Note that  $C_{e_{0}}^{\Delta}=\emptyset$ (see \Cref{claim:BE-properties}(4)) and $N^{2}(e_{0})=T_{4}(e_{0})\cup T_{6}(e_{0})$ (see \Cref{claim:BE-type5-empty}).
Since  $h_{1}\notin N^{2-}(g_{2})$ (see \Cref{claim:h1h2}), we must have 
$s_1v,s_1v_1,t_1v,t_1v_1\notin E(G)$ and so 
$f_{1}\in T_{6}(e_{0})$.
Similarly, since  $h_{2}\notin N^{2-}(g_{1})$ (see \Cref{claim:h1h2}), we  have 
$s_2u,s_2u_1,t_2u,t_2u_1\notin E(G)$ and so 
$f_{2}\in T_{6}(e_{0})$.
 Recall that $|N(g_{2})\cap T_{6}(e)|\ge2$, there exists one edge $f_{3}=v_{1}s_{3}\in (N(g_{2})\cap T_{6}(e_{0}))\setminus\{f_{2}\}$.
 It is possible that $e_{2}\in \{f_{2},f_{3}\}$.
 However, this will not affect the following arguments.

For convenience,  let  $\alpha_{i}=\phi(f_{i})$ for each $i\in[1,3]$.
Note that  $\{\alpha_{0},\alpha_{1},\alpha_{2},\alpha_{3}\}\subseteq \phi(T_6(e_0))$ and  $\{\beta_{1},\beta_{2}\}\subseteq \phi(N(e_0))$. Hence,
by \Cref{claim:BE-properties}, it holds that  $\{\alpha_{0},\alpha_{1},\alpha_{2},\alpha_{3}\}\cap \{\beta_{1},\beta_{2}\}=\emptyset$, $\beta_{1}\neq \beta_{2}$ and $\alpha_{2}\neq\alpha_{3}$.
Moreover, it follows from   $f_{2}\in T_{6}(e_{0})$ that  $Q_{3}=uvv_{1}s_{2}$ is an induced path of length $3$ in $G$ with $uv$ being a bad edge, where $g_{2}=vv_{1}$ and $f_{2}=v_{1}s_{2}$.
     Since  $\phi(g_{2})=\beta_{2}$, by \Cref{claim:path-2}(3), there is exactly one edge $p_{1}=r_{1}r_{1}'$ in $N^{2}_{s_{2}}(f_{2})$ being colored $\beta_{2}$  and $p_{1}\in T_{6}(f_{2})$,  where $r_{1}\in N(s_{2})$. 
    Analogously, as $f_{3}\in T_{6}(e_{0})$ and $Q_{4}=uvv_{1}s_{3}$ is an induced path in $G$,   exactly one edge $p_{2}=r_{2}r_{2}'$ in $N^{2}_{s_{3}}(f_{3})$ with  $r_{2}\in N(s_{3})$ is colored $\beta_{2}$ and $p_{2}\in T_{6}(f_{3})$.   
    Recall that  $f_{2}\in T_{6}(e_{0}) \cap N_{v}^{2}(e_{0})$ and $\phi(f_{2})=\alpha_{2}$,  
    by Observation \ref{Obser:bad edge-type6}, there is exactly one edge $f_{2}'=u'x'$ in $T_{6}(e_{0}) \cap N^{2}_{u}(e_{0})$ with the color $\alpha_{2}$ under $\phi$, where $u'\in N(u)$.
 Possibly   $f_{2}'=f_{1}$. But whether they   are distinct or not will not affect the following arguments.
    We then prove the claim below.

   \begin{claim} \label{claim:v1w}
   	Let $w$ be any vertex in $N(v_{1})$. 
   	If $v_{1}w\in T_{6}(e_{0})$, then $s_{1}w\notin E(G)$.
   \end{claim}
   \begin{proof}
   	Let  $w$ be a vertex in $N(v_{1})$ such that  $v_{1}w\in T_{6}(e_{0})$. 
Assume that $s_{1}w\in E(G)$. 
 Since $f_{1}=u_{1}s_{1}\in T_{6}(e_{0})$,
     $Q_{5}=vuu_{1}s_{1}$ is an induced path  of length $3$ with $vu$ being a bad edge.
      By applying \Cref{claim:path-2}(1) on  $Q_{5}=vuu_{1}s_{1}$, we obtain $C_{f_{1}}^{\Delta}=\emptyset$. 
      Thus, we have $u_{1}w\notin E(G)$ as $s_{1}w\in E(G)$.
      Moreover, 
it follows from $v_{1}w\in T_{6}(e_{0})$  that $uw\notin E(G)$ and $vw\notin E(G)$.    
      Therefore, $Q_{6}=vuu_{1}s_{1}w$ is also an induced path in $G$. 
      Since $h_{1}\notin N^{2-}(g_{2})$ (see \Cref{claim:h1h2}),  we have $h_{1}\neq s_{1}w$ and so $h_{1}\in N(s_{1}w)\setminus\{f_{1}\}$. 
      As  $g_{2}\in N_{w}^{2}(s_{1}w)$,
then   by \Cref{claim:path-1},  we must have  $\phi(h_{1})\neq\phi(g_{2})$. This contradicts the fact that  $\phi(h_{1})=\phi(g_{2})=\beta_{2}$. 
        Hence,  $s_{1}w\notin E(G)$.
    \end{proof}
   
   \Cref{claim:v1w} implies immediately that $s_{1}s_{2}, s_{1}s_{3}\notin E(G)$.
    The remainder of  the proof is divided into the following two cases according to whether the edges $p_{1}$ and   $p_{2}$   are different from $h_{1}$.

  {\bf Case 1.} 	$p_{1}\neq h_{1}$ or $p_{2}\neq h_{1}$.     
  
  By symmetry, we may assume that $p_{1}\neq h_{1}$.
It follows that $s_{2}t_{1}\notin E(G)$.
Recall that $s_1v_1,t_1v_1\notin E(G)$ and $s_{1}s_{2}\notin E(G)$, and thus the distance between $h_{1}$ and $f_{2}$ is  at least $3$ in $G$.
We can obtain a new coloring  $\sigma$ of $G$ by recoloring $g_{1}$ and $f_{2}$ with the same color  $\beta_{2}$, $g_{2}$ with $\beta_{1}$ and  $e_{0}$ with $\alpha_{2}$ (possibly $\alpha_{2}=\alpha_{0}$),  as illustrated in \Cref{fig:sigma1}. 
   In the following, we will show that  $\sigma$ is a good coloring of $G$ with $\kappa_{1}(\sigma)< \kappa_{1}(\phi)$, which contradicts the $1$-optimality of $\phi$.
   
      	\begin{figure}[htbp]  
   	\centering
   	\resizebox{16cm}{5.5cm}{\includegraphics{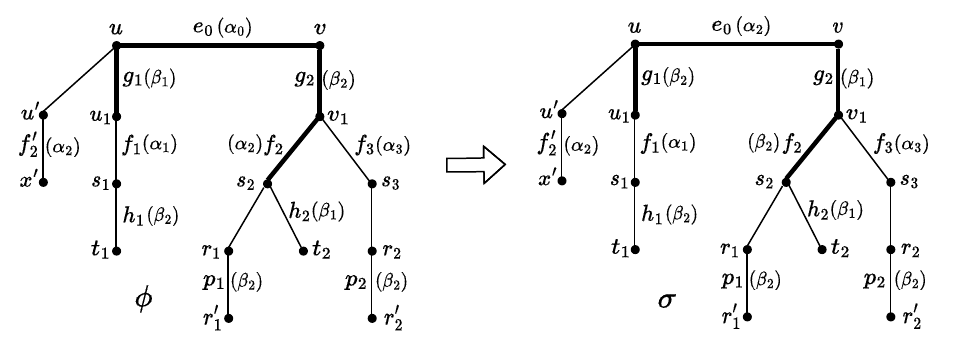}}
   	\caption{The illustration of Case 1}
   	\label{fig:sigma1}
   \end{figure}

 Firstly, note  that   the edge $f_{2}'$ is the only $2$-neighbor of $e_{0}$ colored $\alpha_{2}$ under $\sigma $ and $f_{2}'\in T_{6}(e_{0})$.
  Thus,    the  edge $e_{0}$  is not a bad edge with respect to  $\sigma $. 
  Moreover, since $f_{2}'\in T_{6}(e_{0})$,  $Q_{7}=vuu'x'$ is an induced path of length $3$ in $G$. Then by \Cref{claim:path-2}(2), 
  we have $\phi(f_2')\notin\phi(N^{2-}(f_2')\setminus\{e_0\})$ and so
  the edge $e_{0}$ is the only $2^-$-neighbor of $f_{2}'$ colored $\alpha_{2}$ under $\sigma $. Hence, the  edge $f_{2}'$  is  also not a bad edge with respect to  $\sigma $. 
   
Secondly, since $g_{1},g_{2}\in N(e_{0})$, by \Cref{claim:BE-2-neighborhood}, there is no edge in $N^{2}(e_{0})$ being colored $\beta_{1}$ or $\beta_{2}$ under $\phi$.
Since $u_1v_1,u_1s_2\notin E(G)$, $f_{2}\notin N^{2-}(g_{1})$.
    Recall that $h_{1}=s_{1}t_{1}$ is the only edge in $N^{2}_{u_{1}}(g_{1})$  with $\phi(h_{1})=\beta_{2}$ (as established by applying \Cref{claim:path-1} to $Q_{1}=vuu_{1}$ in the third paragraph of the proof). Hence,   
   $h_{1}$ is the only $2^-$-neighbor of $g_{1}$ colored $\beta_{2}$ under $\sigma$ and $h_1\in T_6(g_1)$ (see \Cref{claim:h1h2}). 
Recall that   $Q_{5}=vuu_{1}s_{1}$ is an induced path. This, together with  	$h_{1}\in T_{6}(g_{1})$ and $h_{1}\notin N^{2-}(e_{0})$ (see \Cref{claim:h1h2}), implies that
$Q_{8}=vuu_{1}s_{1}t_{1}$ is also an induced path in $G$.
     Thus by \Cref{claim:path-2}(2),  we have
     $\phi(h_{1})=\beta_{2}\notin \phi(N^{2-}(h_{1}))$ and so $g_{1}$  is the only $2^-$-neighbor of $h_{1}$ colored $\beta_{2}$ under $\sigma$. 
     Therefore,  both $g_{1}$ and $h_{1}$ are  not  bad edges with respect to  $\sigma$. 
     A similar argument shows that 
     $g_{2}$  and $h_{2}$  are  not  bad edges with respect to  $\sigma$.

 Thirdly, recall that    $Q_{3}=uvv_{1}s_{2}$ is an induced path in $G$ and $g_2=vv_1$,      
  hence by \Cref{claim:path-2}(2), we have $\phi(g_{2})=\beta_{2}\notin\phi(N^{2-}(g_{2}))$. 
  Thus, we must have $p_{1}\notin N^{2-}(g_{2})$ as  $\phi(p_{1})=\beta_{2}$.
   This, together with $p_{1}=r_{1}r_{1}'\in T_{6}(f_{2})$ and $g_2$ is the unique edge in $N^{2-}(e_{0})$ with the color $\beta_2$ under $\phi$,
    implies that $Q_{9}=uvv_{1}s_{2}r_{1}r_{1}'$ is an induced path in $G$. 
Again by  \Cref{claim:path-2}(2), we have  $\phi(p_{1})=\beta_{2}\notin\phi(N^{2-}(p_{1}))$ and so   $f_{2}$   is  the only $2$-neighbor of $p_{1}$ colored $\beta_{2}$ under $\sigma$. 
   Moreover, according to  \Cref{claim:path-1}, it is easy to see that  $p_{1}$ is the only $2$-neighbor of $f_{2}$ colored $\beta_{2}$ under $\sigma$.
Hence, the two edges $f_{2}$ and $p_{1}$ are not bad edges   with respect to $\sigma$.

Finally,
if $p_{2}\in \{h_{1},p_{1}\}$, then $p_{2}$ is obviously not a bad edge   with respect to $\sigma$.  If $p_{2}\notin \{h_{1},p_{1}\}$, then according to a similar argument applied to $Q_{9}=uvv_{1}s_{2}r_{1}r_{1}'$,  the path
$Q_{10}=uvv_{1}s_{3}r_{2}r_{2}'$ is also an induced path in $G$ (note that $p_{2}=r_{2}r_{2}'$).
By 
 \Cref{claim:path-2}(2),  we have $\phi(p_{2})=\beta_{2}\notin\phi(N^{2-}(p_{2}))$  and so
   $p_{2}$ has no $2^{-}$-neighbor with the color $\beta_{2}$ under $\sigma$.
 Thus,   $p_{2}$
   is not  a bad edge   with respect to $\sigma$.

 Hence, the above discussion shows that  
 $\sigma$ is a good coloring of $G$ satisfying $\kappa_{1}(\sigma )< \kappa_{1}(\phi)$.

    {\bf Case 2.} 	$p_{1}=p_{2}= h_{1}$. 
    
    In this case,  we must have $s_{2}t_{1},s_{3}t_{1}\in E(G)$ as $s_{1}s_{2}, s_{1}s_{3}\notin E(G)$     (refer to \Cref{fig:sigma2}). 
     Let $\gamma=\phi(s_{3}t_{1})$.
     Recall that $Q_{2}=uvv_{1}$ is induced,
    we must have $\gamma\notin\{\beta_{1},\beta_{2},\alpha_{2}\}$ due to \Cref{claim:path-1}.
     Now,
     we recolor $g_{1}$  with  $\beta_{2}$, $g_{2}$ with $\beta_{1}$, 
       $e_{0}$  and  $s_{3}t_{1}$ with the same color $\alpha_{2}$ (possibly $\alpha_{2}=\alpha_{0}$) and $f_{2}$ with $\gamma$.   This yields  a new coloring of $G$  called $\sigma$.
       Similar to the arguments in the proof of  Case 1, it
        is easy to check that  $\sigma$ is a good coloring of $G$
        and the six edges $e_{0},f_{2}',g_{1},g_{2},h_{1},h_{2}$ are not bad edges with respect to $\sigma$. 
            \begin{figure}[htbp]  
        	\centering
        	\resizebox{16cm}{5.5cm}{\includegraphics{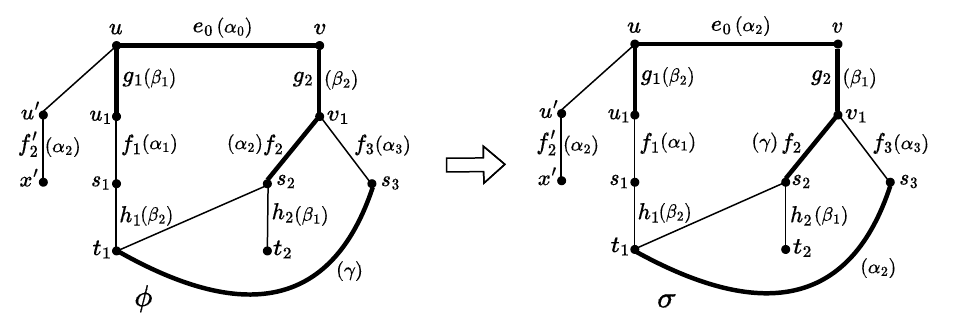}}
        	\caption{The illustration of Case 2}
        	\label{fig:sigma2}
        \end{figure}
     
       Notice that $s_{3}t_{1}\in T_{4}(f_{2})$,  
       and that $Q_{3}=uvv_{1}s_{2}$ is an induced path in $G$. Hence, by \Cref{claim:path-1}, it is easy to check that
       no edge in $N^{2-}(f_{2})\setminus\{s_{3}t_{1},e_{0}\}$ is colored $\gamma$ under $\phi$.
It follows that $f_{2}$ has no  $2^{-}$-neighbor being colored   $\gamma$ under $\sigma$.
   Hence,   $f_{2}$ 
       is not a bad edge with respect to $\sigma$.     
         Recall that $h_{1}=s_{1}t_{1}\notin N^{2-}(e_{0})\cup N^{2-}(g_{2})$ (see \Cref{claim:h1h2}), we have  $t_{1}u,t_{1}v,t_{1}v_{1}\notin E(G)$. 
       This, together with  the induced path $Q_{4}=uvv_{1}s_{3}$, implies that   $Q_{11}=uvv_{1}s_{3}t_{1}$ is  also induced and $s_{3}t_{1}\notin N^{2-}(e_{0})$.
      Then by  \Cref{claim:path-1},
       $f_{2}$ is the only edge in $N^{2-}(s_{3}t_{1})$ being colored $\alpha_{2}$ under $\phi$.
        Therefore, $s_{3}t_{1}$  has no  $2^{-}$-neighbor being colored   $\alpha_{2}$ under $\sigma$ and so it
        is not a bad edge with respect to $\sigma$.

   Thus,   $\sigma$ is a good coloring of $G$ with $\kappa_{1}(\sigma)< \kappa_{1}(\phi)$, contradicting the $1$-optimality of $\phi$ again.

\medskip
We have deduced contradictions in  both cases.
Therefore, there is no bad edge with respect to $\phi$ in $G$.
By Observation \ref{Obser:no bad edge}, 
 $\phi$ is both  a semistrong edge coloring and 
 a $(0,1)$-relaxed strong edge coloring using at most $\Delta^{2}-1$ colors. This completes the proof of this lemma.
\end{proof}

 Theorem \ref{lemma:delta3+} follows from  Lemmas \ref{lemma:G in Gp} and \ref{lemma:G notin Gp}. 
\section{Summary} \label{sec:5}
In this paper,  we proved that the semistrong chromatic index of   a connected graph  with maximum degree $\Delta$ is at most  $\Delta^{2}-1$, except  $C_{7}$ and $K_{\Delta,\Delta}$.
 This upper bound  is  tight, since the upper bound $3$ is the best possible when $\Delta=2$.
Moreover, as indicated by   Lu{\v{z}}ar, Mockov{\v{c}}iakov{\'a} and  Sot{\'a}k  in \cite{LMS2022}, the $5$-prism (see \Cref{fig:5prism}) shows the sharpness of the bound 8 for the case $\Delta=3$. 
However, they did not find infinitely many graphs attaining the bound 8.
Likewise, we have not found any graphs with maximum degree $\Delta\ge4$ whose semistrong chromatic indices are equal to  $\Delta^{2}-1$.

 For  $\Delta=4$, the graph “$C_{7}$-blowup”  constructed as follows has the semistrong chromatic index $\Delta^{2}-2$: the vertex set $V=\cup_{i=0}^{6}V_{i}$ where every $V_{i}$ is an independent set with two vertices;  for any two different integers $i,j\in\{0,1,\dots,6\}$,   $V_{i}\cap V_{j}=\emptyset$; and every vertex in $V_{i}$ is adjacent to every vertex in $V_{i+1}$, with indices taken modulo 7. 
Therefore, we  believe that  the upper bound $\Delta^{2}-1$ can be further improved.
After some exploration, we propose the following problem.

\medskip
\noindent {\bf Problem 1}: Let $G$ be a connected graph with maximum degree $\Delta$ that is not isomorphic to  $K_{\Delta,\Delta}$.  If  $\Delta$ is  appropriately large, is it true that  $\chi'_{ss}(G)\le \Delta^{2}-\Delta+1$?

\medskip
It should be pointed out that, the above upper bound  if proven, would be the best possible.
Let $H$ denote the graph obtained by taking two copies of  $K_{\Delta-1,\Delta}$ and  adding one edge between two distinct vertices of degree $\Delta-1$ from each of the two copies.
Clearly, the maximum degree of   $H$ is $\Delta$.
Moreover, it is easy to check that $\chi'_{ss}(H)= \Delta^{2}-\Delta+1$.
It follows that,
any graph $G$ with maximum degree $\Delta$ containing $H$ as a subgraph has the semistrong chromatic index at least $\Delta^{2}-\Delta+1$.

Meanwhile, we also proved that 
any  connected graph  with maximum degree $\Delta$, except  $C_{7}$, has 
 $(0,1)$-relaxed  strong chromatic index  at most  $\Delta^{2}-1$.
However, we tried without success finding a graph whose  $(0,1)$-relaxed  strong chromatic index is close to $\Delta^{2}-1$.
We therefore strongly believe that this upper bound is not tight and propose the following conjecture.

\begin{conjecture}\label{Conj-ending}
	For every connected graph $G$ with maximum degree $\Delta$ other than $C_{7}$, 
	\begin{equation*}
		\chi'_{(0,1)}(G)\le\begin{cases}
			\begin{array}{cl}
			\lceil 	\dfrac{5}{8}\Delta^{2}\rceil,                                 & \text{if}\  \Delta \ \text{is even,} \\
				\lceil 	\dfrac{5}{8}\Delta^{2}-\dfrac{1}{4}\Delta+\dfrac{1}{8}\rceil, & \text{if}\  \Delta\  \text{is odd.}
			\end{array}
		\end{cases}	
	\end{equation*}
\end{conjecture}

The graphs  “$C_{5}$-blowups”  constructed by 
Erd\H{o}s and Ne\v{s}et\v{r}il \cite{E1988,EN1989}  indicate that the bounds given in Conjecture \ref{Conj-ending}, if proven, would be tight.
 In some sense, Conjecture \ref{Conj-ending} can be seen as a reinforcement of Erd\H{o}s and Ne\v{s}et\v{r}il's conjecture (see Conjecture \ref{Conj-EN}).
Given a $(0,1)$-relaxed strong $k$-edge-coloring of $G$,   a strong $2k$-edge-coloring of $G$ can be easily obtained  by dividing each color class into two.  This reveals that
if Conjecture \ref{Conj-ending} is true,  then the strong chromatic index of any graph does not exceed the upper bound conjectured by Erd\H{o}s and Ne\v{s}et\v{r}il plus one.  This yields a new way to break though  Erd\H{o}s and Ne\v{s}et\v{r}il's conjecture.\\

 \section*{Acknowledgments}
 The first author was supported by the China Scholarship Council (CSC)	and  SEU Innovation Capability Enhancement Plan for Doctoral Students (CXJH\_SEU 24119).
The second author was supported by National Natural Science Foundation of China 11771080.
We sincerely thank the anonymous reviewers for their careful reading and insightful comments, which have significantly improved this work.

 \section*{Data Availability Statement}
Data sharing not applicable to this article as no datasets were generated or analyzed during the current study.

\bibliographystyle{amsplain}
\bibliography{ref}

\end{document}